\documentclass{SCAE}
\numberwithin{equation}{section}

\newcommand{\nn}{\nonumber}

\newcommand{\vep}{\varepsilon}

\newcommand{\bx}{{\bf x} }

\newcommand{\bB}{{\bf B} }
\newcommand{\bk}{{\bf k} }
\newcommand{\bJ}{{\bf J} }
\newcommand{\bee}{{\bf e} }
\newcommand{\beU}{{ U} }
\newcommand{\bV}{{\bf V}}

\newcommand{\bG}{{\bf G} }
\newcommand{\bW}{{\bf W} }
\newcommand{\bF}{{\bf F} }

\newcommand{\eps}{\varepsilon}

\newcommand{\be}{\begin{equation}}
\newcommand{\ee}{\end{equation}}
\newcommand{\ba}{\begin{array}}
\newcommand{\ea}{\end{array}}
\newcommand{\bea}{\begin{eqnarray}}
\newcommand{\eea}{\end{eqnarray}}
\newcommand{\beas}{\begin{eqnarray*}}
\newcommand{\eeas}{\end{eqnarray*}}

\begin{document}

\Year{2016} %
\Month{January}
\Vol{59} %
\No{1} %
\BeginPage{1} %
\EndPage{XX} %
\AuthorMark{Bao W {\it et al.}}
\ReceivedDay{XXXX}
\AcceptedDay{XXXX}
\PublishedOnlineDay{XXXX}
\DOI{10.1007/s11425-000-0000-0} 

\title[Error estimates for nonlinear Dirac equation]{Error estimates of numerical methods for \\
the nonlinear Dirac equation \\
in the nonrelativistic limit regime}{}


\author[1]{BAO WeiZhu}{Corresponding author. {\rm http://www.math.nus.edu.sg/\~{}bao/}}
\author[2,3]{CAI YongYong}{}
\author[1]{JIA XiaoWei}{}
\author[4]{YIN Jia}{}

\address[{\rm1}]{Department of Mathematics,  National University of
Singapore, Singapore {\rm 119076}, Singapore}
\address[{\rm2}]{Beijing Computational Science Research Center,
No. 10 West Dongbeiwang Road, Beijing {\rm 100193}, P. R. China}
\address[{\rm3}]{Department of Mathematics, Purdue University, West Lafayette,
IN {\rm 47907}, USA}
\address[{\rm4}]{NUS Graduate School for Integrative Sciences and Engineering (NGS), National University of
Singapore, Singapore {\rm 117456}, Singapore}
\Emails{ matbaowz@nus.edu.sg,
yongyong.cai@gmail.com, A0068124@nus.edu.sg, e0005518@u.nus.edu}\maketitle


 {\begin{center}
\parbox{14.5cm}{\begin{abstract}
We present several numerical methods and establish their error estimates for the
discretization of the nonlinear Dirac equation
in the nonrelativistic limit regime, involving a small dimensionless parameter $0<\varepsilon\ll 1$
which is inversely proportional to the speed of light. In this limit
regime, the solution is highly oscillatory in time, i.e. there are propagating waves
with wavelength $O(\varepsilon^2)$ and $O(1)$ in time and space, respectively.
We begin with the conservative Crank-Nicolson finite difference (CNFD)
method and establish rigorously its error estimate
which depends explicitly on the mesh size $h$ and time step $\tau$ as well
as the small parameter $0<\varepsilon\le 1$.
Based on the error bound, in order to obtain `correct'
numerical solutions in the nonrelativistic limit regime, i.e. $0<\varepsilon\ll 1$,
the CNFD method requests the $\varepsilon$-scalability:
$\tau=O(\varepsilon^3)$ and $h=O(\sqrt{\varepsilon})$.
Then we propose and analyze two numerical methods
for the discretization of the nonlinear Dirac equation by using the Fourier
spectral discretization for spatial derivatives combined with the
exponential wave integrator and time-splitting technique
for temporal derivatives, respectively.
Rigorous error bounds for the two numerical methods
show that their $\varepsilon$-scalability is improved
to $\tau=O(\varepsilon^2)$ and $h=O(1)$ when $0<\varepsilon\ll 1$
compared with the CNFD method.
Extensive numerical results
are reported to confirm our error estimates.
\end{abstract}}\end{center}}



\renewcommand{\baselinestretch}{1.2}
\begin{center} \renewcommand{\arraystretch}{1.5}
{\begin{tabular}{lp{0.8\textwidth}} \hline \scriptsize
{\bf Citation:}\!\!\!\!&\scriptsize Bao W, Cai Y, Jia X, Yin J. \makeatletter\@titlehead.
Sci China Math, 2016, 59,
 doi:~\@DOI\makeatother\vspace{1mm}
\\
\hline
\end{tabular}}\end{center}

\baselineskip 11pt\parindent=10.8pt  \wuhao

\section{Introduction}\setcounter{equation}{0}
In particle physics and/or relativistic quantum mechanics,
the Dirac equation, which was derived
by the British physicist Paul Dirac in 1928 \cite{Dirac1,Dirac3},
is a relativistic wave equation for describing all
spin-$1/2$ massive particles, such as electrons and positrons.
It is consistent with both the principle of quantum
mechanics and the theory of special relativity, and
was the first theory to fully account for relativity in the context
of quantum mechanics. It accounted for the fine details
of the hydrogen spectrum in a completely rigorous way
and provided the entailed explanation of spin as a consequence
of the union of quantum mechanics and relativity -- and
the eventual discovery of the positron -- represent one of
the greatest triumphs of theoretical physics.
Since the graphene was first produced in the lab in 2003
\cite{AMPGMKWTNLG,NGPNG}, the Dirac equation has been extensively adopted to
study theoretically the structures and/or dynamical
properties of graphene and graphite as well as other two dimensional
(2D) materials \cite{AMPGMKWTNLG,FW,NGPNG}. Recently, the Dirac equation
has also been adopted to study the
relativistic effects in molecules in super intense lasers, e.g. attosecond lasers \cite{FLB2,FLB3}
and the motion of nucleons in the covariant density function theory
in the relativistic framework \cite{Ring,Meng}.

The nonlinear Dirac equation (NLDE), which was proposed in 1938 \cite{Ivan},
is a model of self-interacting Dirac fermions
in quantum field theory \cite{FLR,FS,Soler,Thirring} and is widely considered
as a toy model of self-interacting spinor fields in quantum physics \cite{Soler,Thirring}.
In fact, it also appears in the
Einstein-Cartan-Sciama-Kibble theory of gravity, which
extends general relativity to matter with intrinsic angular
momentum (spin) \cite{Heis}.
In the resulting field equations, the minimal coupling
between the homogeneous torsion tensor and Dirac spinors
generates an axial-axial, spin-spin interaction in fermionic matter,
which becomes nonlinear (cubic) in the spinor field and
significant only at extremely high densities.
Recently, the NLDE has been adapted
as a mean field model for Bose-Einstein condensates (BECs) \cite{CEL,HC,HWC}
and/or cosmology \cite{Sah}.
In fact, the experimental advances in  BECs and/or graphene as well as 2D materials
have renewed extensively
the research interests on the mathematical analysis and numerical simulations
of the Dirac equation and/or the NLDE without/with
electromagnetic potentials, especially the honeycomb lattice potential \cite{AZ,FW,FW2}.

Consider the NLDE in three dimensions (3D) for
describing the dynamics of spin-$1/2$
self-interacting massive Dirac fermions within
external time-dependent electromagnetic potentials \cite{Dirac1,FLR,FS,Soler,Thirring,HC,HWC}
\begin{equation}
\label{3DME}
i\hbar\partial_t\Psi = \Bigl[-ic\hbar\sum_{j=1}^3\alpha_j\partial_j
+mc^2\beta \Bigr] \Psi+e\Bigl[V(t,\bold{x})I_4 -\sum_{j=1}^3A_j(t,\mathbf{x})\alpha_j
\Bigr]\Psi+\bF(\Psi)\Psi,\quad \bx\in{\mathbb R}^3,
\end{equation}
where $i=\sqrt{-1}$,  $t$ is time, $\bx=(x_1,x_2,x_3)^T\in {\mathbb R}^3$ (equivalently written as $\bx=(x,y,z)^T$)
is the spatial coordinate vector, $\partial_k=\frac{\partial}{\partial x_k}$ ($k=1,2,3$),
 $\Psi :=\Psi(t,\bx)=(\psi_1(t,\bx),\psi_2(t,\bx), \psi_3(t,\bx), \psi_4(t,\bx))^T\in\mathbb{C}^4$
 is the complex-valued vector
wave function of the ``spinorfield''.
$I_n$ is the $n\times n$ identity matrix for $n\in {\mathbb N}$,
$V:=V(t,\bx)$ is the real-valued electrical potential and
${\bf A}:={\bf A}(t,\bx)=(A_1(t,\bx), A_2(t,\bx), A_3(t,\bx))^T$ is the real-valued magnetic potential vector,
and hence the electric field is given by ${\bf E}(t,\bx)=-\nabla V-\partial_t {\bf A}$ and
the magnetic field is given by ${\bf B}(t,\bx)={\rm curl}\, {\bf A}=\nabla \times {\bf A}$.
Different cubic nonlinearities have been proposed
and studied for the NLDE (\ref{3DME}) from different applications \cite{FLR,FS,Soler,Thirring,HC,HWC}.
For the simplicity of notations, here we take
$\bF(\Psi)=g_1 \left(\Psi^*\beta\Psi\right)\beta+g_2|\Psi|^2I_4$ with $g_1,g_2\in{\mathbb R}$ two constants
and $\Psi^*=\overline{\Psi}^T$, while $\overline{f}$ denotes the complex conjugate of  $f$,
which is motivated  from the so-called Soler model, e.g. $g_2=0$ and $g_1\ne0$,
in quantum field theory \cite{FLR,FS,Soler,Thirring} and BECs with a chiral confinement and/or spin-orbit coupling,
e.g. $g_1=0$ and $g_2\ne0$ \cite{CEL,HC,HWC}. We remark here that our numerical methods and their error estimates
can be easily extended to the NLDE with other nonlinearities \cite{Thirring,Sah,Shao}.
The physical constants are: $c$ for the speed of light, $m$ for the particle's rest mass,
$\hbar$ for the Planck constant and $e$ for the unit charge. In addition, the $4\times 4$
matrices $\alpha_1$, $\alpha_2$, $\alpha_3$ and $\beta$ are defined as
\be \label{alpha}
\alpha_1=\left(\begin{array}{cc}
\mathbf{0} & \sigma_1  \\
\sigma_1 & \mathbf{0}  \\
\end{array}
\right),\qquad
\alpha_2=\left(\begin{array}{cc}
\mathbf{0} & \sigma_2 \\
\sigma_2 & \mathbf{0} \\
\end{array}
\right), \qquad
\alpha_3=\left(\begin{array}{cc}
\mathbf{0} & \sigma_3 \\
\sigma_3 & \mathbf{0} \\
\end{array}
\right),\qquad
\beta=\left(\begin{array}{cc}
I_{2}& \mathbf{0} \\
\mathbf{0} & -I_{2} \\
\end{array}
\right),
\ee
with  $\sigma_1$, $\sigma_2$,  $\sigma_3$
(equivalently written $\sigma_x$, $\sigma_y$, $\sigma_z$) being the Pauli matrices  defined as
\be\label{Paulim}
\sigma_{1}=\left(
\begin{array}{cc}
0 & 1  \\
1 & 0  \\
\end{array}
\right), \qquad
\sigma_{2}=\left(
\begin{array}{cc}
0 & -i \\
i & 0 \\
\end{array}
\right),\qquad
\sigma_{3}=\left(
\begin{array}{cc}
1 & 0 \\
0 & -1 \\
\end{array}
\right).
\ee

Similar to the Dirac equation \cite{BCJ}, by a proper
nondimensionalization (with the choice of $x_s$, $t_s=\frac{mx_s^2}{\hbar}$,
$A_s=\frac{m v^2}{e}$ and $\psi_s=x_s^{-3/2}$ as
the dimensionless length unit, time unit, potential unit and spinor field unit, respectively)
and dimension reduction \cite{BCJ}, we can obtain the
dimensionless NLDE in $d$-dimensions ($d=3,2,1$)
\be
\label{SDEd}
i\partial_t\Psi=\Bigl[-\frac{i}{\varepsilon}\sum_{j=1}^{d}\alpha_j
\partial_j+\frac{1}{\varepsilon^2}\beta\Bigr]\Psi
+\Bigl[V(t,\bx)I_4-\sum_{j=1}^{d}A_j(t,\bx)\alpha_j\Bigr]\Psi+\bF(\Psi)\Psi,
\quad \bx\in{\mathbb R}^d,
\ee
where $\varepsilon$ is a dimensionless parameter  inversely proportional to the speed of light given by
\begin{equation}\label{eps}
0<\varepsilon := \frac{x_s}{t_s\,c}=\frac{v}{c}\le 1,
\end{equation}
with $v=\frac{x_s}{t_s}$ the wave speed, and
\be\label{nolc}
\bF(\Psi)=\lambda_1 \left(\Psi^*\beta\Psi\right)\beta+\lambda_2|\Psi|^2I_4,\qquad \Psi\in {\mathbb C}^4,
\ee
with  $\lambda_1=\frac{g_1}{mv^2x_s^3}\in{\mathbb R}$ and $\lambda_2=\frac{g_2}{mv^2x_s^3}\in{\mathbb R}$
two dimensionless constants for the interaction strength.

For the dynamics, the initial condition is given as
\begin{equation*}
\Psi(t=0,\bx)=\Psi_0(\bx), \qquad \bx\in\mathbb R^d.
\end{equation*}

The NLDE (\ref{SDEd}) is dispersive and time symmetric \cite{XST}.
Introducing the position density $\rho_j$ for the $j$-component ($j=1,2,3,4$), 
the total density $\rho$ as well as the current density
$\bJ(t,\bx)=(J_1(t,\bx),J_2(t,\bx)$, $J_3(t,\bx))^T$
\begin{equation} \label{obser11}
\rho(t,\bx)=\sum_{j=1}^4\rho_j(t,\bx)=\Psi^*\Psi,
\ \ \rho_j(t,\bx)=|\psi_j(t,\bx)|^2, \quad 1\leq j\leq 4;\ \
J_l(t,\bx)=\frac{1}{\varepsilon}\Psi^*\alpha_l\Psi,\quad l=1,2,3,
\end{equation}
then the following conservation law can be obtained from the NLDE
(\ref{SDEd})
\begin{equation} \label{cons}
\partial_t\rho(t,\bx)+\nabla\cdot \bJ(t,\bx) = 0, \qquad \bx\in{\mathbb R}^d, \quad t>0.
\end{equation}
Thus the NLDE (\ref{SDEd}) conserves the total mass as
\be\label{norm}
\|\Psi(t,\cdot)\|^2:=\int_{{\mathbb R}^d}|\Psi(t,\bx)|^2\,d\bx=
\int_{{\mathbb R}^d}\sum_{j=1}^4|\psi_j(t,\bx)|^2\,d\bx\equiv \|\Psi(0,\cdot)\|^2
=\|\Psi_0\|^2, \qquad t\ge0.
\ee
If the electric potential $V$ is perturbed by a constant, e.g. $V(t,\bx)\to V(t,\bx)+V^0$ with
$V^0$ being a real constant, then the solution $\Psi(t,\bx)\to e^{-iV^0t}\Psi(t,\bx)$  which
implies the density of each component $\rho_j$ ($j=1,2,3,4$) and the total
density $\rho$ unchanged. When $d=1$, if the magnetic potential $A_1$
is perturbed by a constant, e.g. $A_1(t,\bx)\to A_1(t,\bx)+A_1^0$ with
$A_1^0$ being a real constant, then the solution $\Psi(t,\bx)\to e^{iA_1^0t\alpha_1}\Psi(t,\bx)$  which
implies the total density $\rho$ unchanged; but this property is not valid when $d=2,3$.
These properties are usually called as time transverse invariant.
In addition, when the electromagnetic potentials are time-independent,
i.e. $V(t,\bx)=V(\bx)$ and $A_j(t,\bx)=A_j(\bx)$ for $j=1,2,3$,
 the following energy functional is also conserved
\bea\label{engery60}
E(t)&:=&\int_{\mathbb{R}^d}\left[-\frac{i}{\varepsilon}
\sum_{j=1}^d\Psi^*\alpha_j\partial_j\Psi+\frac{1}{\varepsilon^2}\Psi^*\beta\Psi+V(\bx)|\Psi|^2+
G(\Psi)-\sum_{j=1}^dA_j(\bx)\Psi^*\alpha_j\Psi\right]d\bx\nonumber\\
&\equiv&E(0),\qquad t\ge0,
\eea
where
\be\label{GPsi}
G(\Psi)=\frac{\lambda_1}{2} \left(\Psi^*\beta\Psi\right)^2+\frac{\lambda_2}{2}|\Psi|^4,
\qquad \Psi\in {\mathbb C}^4.
\ee
Furthermore,  if the external electromagnetic potentials are constants, i.e.
$V(t,\bx)\equiv V^0$ and $A_j(t,\bx)\equiv A_j^0$ for $j=1,2,3$,
the NLDE (\ref{SDEd}) admits the plane wave solution as
$\Psi(t,\bx)={\bf B}\,e^{i(\bk\cdot\bx-\omega t)}$,
where the time frequency $\omega$, amplitude vector $\bB\in {\mathbb R}^4$
and spatial wave number $\bk=(k_1,\ldots,k_d)^T\in {\mathbb R}^d$ satisfy
\be\label{disp}
\omega \bB=\left[\sum_{j=1}^{d}\left(\frac{k_j}{\varepsilon}-A_j^0\right)\alpha_j
+\frac{1}{\varepsilon^2}\beta
+V^0I_{4}+\lambda_1 \left(\bB^*\beta \bB\right)\beta+\lambda_2|\bB|^2I_4\right]\bB,
\ee
which immediately implies the {\sl dispersion relation} of the NLDE (\ref{SDEd}) as
\be\label{disperg}
\omega:=\omega(\bk,\bB)=V^0+\lambda_2|\bB|^2
\pm \frac{1}{\vep^2}\sqrt{\left[1+\vep^2
\lambda_1 \left(\bB^*\beta \bB\right)\right]^2+\vep^2\left|\bk-\vep {\bf A}^0\right|^2}=O\left(\frac{1}{\vep^2}\right),\quad
\bk\in {\mathbb R}^d.
\ee

Again, similar to the Dirac equation \cite{BCJ}, in several applications
in one dimension (1D) and two dimensions (2D), the NLDE (\ref{SDEd})
can be simplified to the following NLDE in $d$-dimensions ($d=1,2$) with
$\Phi:=\Phi(t,\bx)=(\phi_1(t,\bx),\phi_2(t,\bx))^T \in{\mathbb C}^2$ \cite{FLR,FS,Soler}
\be
\label{SDEdd}
i\partial_t\Phi=\Bigl[-\frac{i}{\varepsilon}\sum_{j=1}^{d}\sigma_j
\partial_j+\frac{1}{\varepsilon^2}\sigma_3\Bigr]\Phi
+\Bigl[V(t,\bx)I_2-\sum_{j=1}^{d}A_j(t,\bx)\sigma_j\Bigr]\Phi+\bF(\Phi)\Phi,
\quad \bx\in{\mathbb R}^d,
\ee
where
\be\label{nolcc}
\bF(\Phi)=\lambda_1 \left(\Phi^*\sigma_3\Phi\right)\sigma_3+\lambda_2|\Phi|^2I_2,\qquad \Phi\in {\mathbb C}^2,
\ee
with $\lambda_1\in{\mathbb R}$ and $\lambda_2\in{\mathbb R}$
two dimensionless constants for the interaction strength.
Again, the initial condition for dynamics is given as
\be\label{SDEini}
\Phi(t=0,\bx)=\Phi_0(\bx), \qquad \bx\in{\mathbb R}^d.
\ee
The NLDE (\ref{SDEdd}) is dispersive and time symmetric.
By introducing the position density $\rho_j$ for the $j$-th component ($j=1,2$), 
the total density $\rho$ as well as the current density $\bJ(t,\bx)=(J_1(t,\bx),J_2(t,\bx))^T$
\be\label{obser12}
\rho(t,\bx)=\sum_{j=1}^2\rho_j(t,\bx)=\Phi^\ast \Phi,
\quad \rho_j(t,\bx)=|\phi_j(t,\bx)|^2, \quad
J_j(t,\bx)=\frac{1}{\varepsilon}\Phi^*\sigma_j\Phi, \quad j=1,2,
\ee
 the conservation law (\ref{cons})
is also satisfied \cite{BHM}.  In addition, the NLDE (\ref{SDEdd})
conserves the total mass as
\be\label{normd}
\|\Phi(t,\cdot)\|^2:=\int_{{\mathbb R}^d}|\Phi(t,\bx)|^2\,d\bx=\int_{{\mathbb R}^d}\sum_{j=1}^2|\phi_j(t,\bx)|^2\,d\bx\equiv \|\Phi(0,\cdot)\|^2=\|\Phi_0\|^2, \qquad t\ge0.
\ee
Again, if the electric potential $V$ is perturbed by a constant, e.g. $V(t,\bx)\to V(t,\bx)+V^0$ with
$V^0$ being a real constant,  the solution $\Phi(t,\bx)\to e^{-iV^0t}\Phi(t,\bx)$  which
implies the density of each component $\rho_j$ ($j=1,2$) and the total
density $\rho$ unchanged. When $d=1$, if the magnetic potential $A_1$
is perturbed by a constant, e.g. $A_1(t,\bx)\to A_1(t,\bx)+A_1^0$ with
$A_1^0$ being a real constant, the solution $\Phi(t,\bx)\to e^{iA_1^0t\sigma_1}\Phi(t,\bx)$
implying the total density $\rho$ unchanged; but this property is not valid when $d=2$.
When the electromagnetic potentials are time-independent,
i.e. $V(t,\bx)=V(\bx)$ and $A_j(t,\bx)=A_j(\bx)$ for $j=1,2$,
 the following energy functional is also conserved
\bea \label{engery65}
E(t)&:=&\int_{\mathbb{R}^d}\left(-\frac{i}{\varepsilon}
\sum_{j=1}^d\Phi^*\sigma_j\partial_j\Phi+\frac{1}{\varepsilon^2}\Phi^*\sigma_3\Phi+V(\bx)|\Phi|^2-
\sum_{j=1}^dA_j(\bx)\Phi^*\sigma_j\Phi+G(\Phi)\right)d\bx\nn\nonumber\\
&\equiv&E(0),\qquad t\ge0,
\eea
where
\be\label{GPhi}
G(\Phi)=\frac{\lambda_1}{2} \left(\Phi^*\sigma_3\Phi\right)^2+\frac{\lambda_2}{2}|\Phi|^4,\qquad \Phi\in{\mathbb C}^2.
\ee
Furthermore,  if the external electromagnetic potentials are constants, i.e.
$V(t,\bx)\equiv V^0$ and $A_j(t,\bx)\equiv A_j^0$ for $j=1,2$,
the NLDE (\ref{SDEdd}) admits the plane wave solution as
$\Phi(t,\bx)={\bf B}\,e^{i(\bk\cdot\bx-\omega t)}$,
where the time frequency $\omega$, amplitude vector $\bB\in {\mathbb R}^2$
and spatial wave number $\bk=(k_1,\ldots,k_d)^T\in {\mathbb R}^d$ satisfy
\be\label{dispp}
\omega \bB=\Bigl[\ \sum_{j=1}^{d}\left(\frac{k_j}{\varepsilon}-A_j^0\right)\sigma_j
+\frac{1}{\varepsilon^2}\sigma_3
+V^0I_{2}+\lambda_1 \left(\bB^*\sigma_3 \bB\right)\sigma_3+\lambda_2|\bB|^2I_2\Bigr]\bB,
\ee
which immediately implies the {\sl dispersion relation} of the NLDE (\ref{SDEdd}) as
\be\label{dispergg}
\omega:=\omega(\bk,\bB)=V^0+\lambda_2|\bB|^2
\pm \frac{1}{\vep^2}\sqrt{\left[1+\vep^2
\lambda_1 \left(\bB^*\sigma_3 \bB\right)\right]^2+\vep^2\left|\bk-\vep {\bf A}^0\right|^2}=
O\left(\frac{1}{\vep^2}\right),\quad
\bk\in {\mathbb R}^d.
\ee

\begin{figure}[htb]
\centerline{\includegraphics[height=5cm,width=13cm]{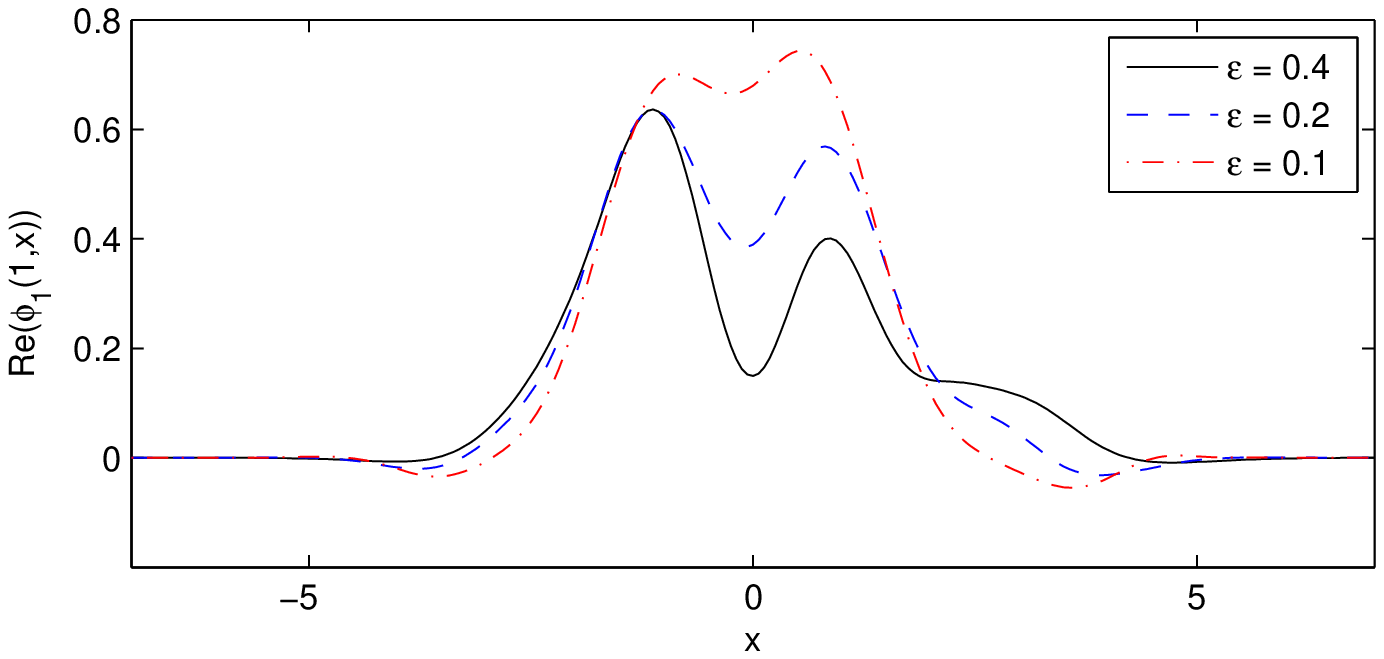}}
\centerline{\includegraphics[,height=5cm,width=13cm]{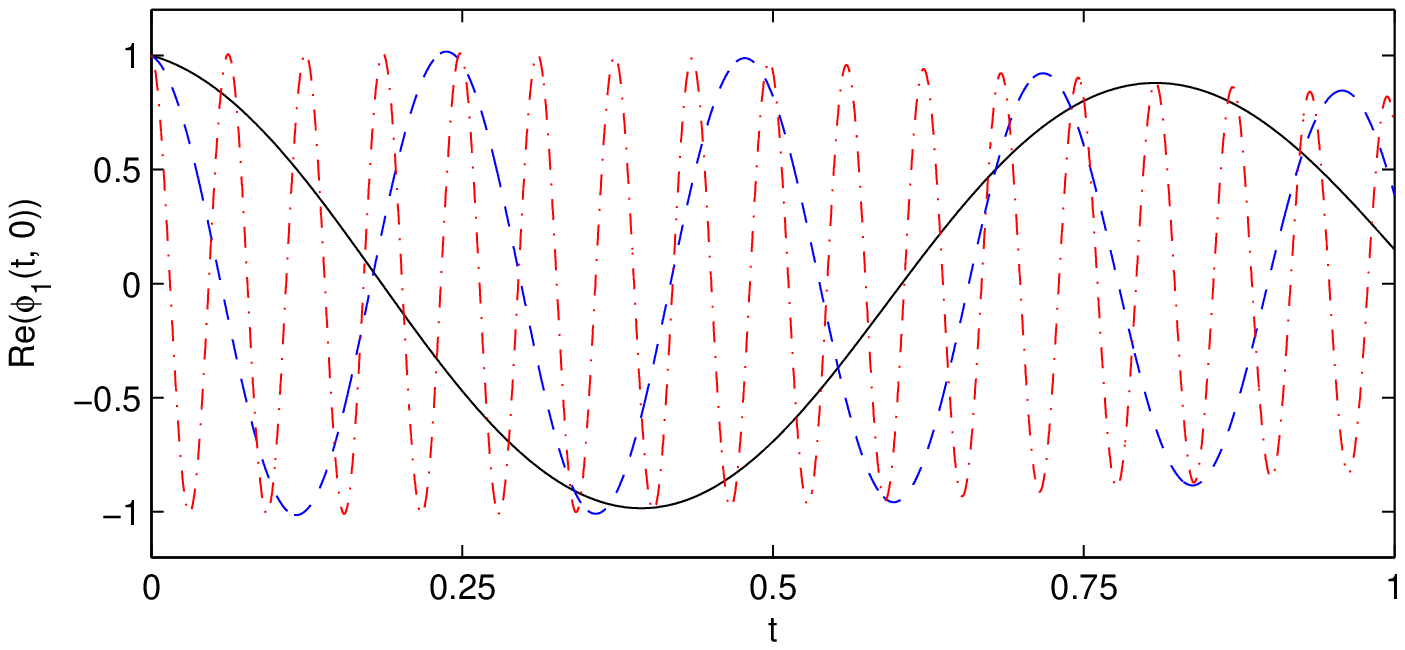}}
\caption{The solution $\phi_1(t=1,x)$ and $\phi_1(t,x=0)$ of the NLDE
(\ref{SDEdd}) with $d=1$ for different $\varepsilon$. ${\rm Re}(f)$ denotes the real part
of $f$.}
\label{Oscillation}
\end{figure}

For the NLDE (\ref{SDEd}) (or (\ref{SDEdd}))
with $\varepsilon=1$, i.e. $O(1)$-speed of light regime,
there are extensive analytical and numerical results
in the literatures. For the existence and multiplicity of bound states and/or standing wave solutions,
we refer to \cite{Bala,Bala1,Bar,Caz,Dol,Esteban0,Esteban,Komech} and references therein.
Particularly, when $d=1$, $\varepsilon=1$, $V(t,x)\equiv0$ and $A_1(t,x)\equiv 0$ in (\ref{SDEdd})
and $\lambda_1=-1$ and $\lambda_2=0$ in (\ref{nolcc}), the NLDE (\ref{SDEdd}) admits
soliton solutions which was given explicitly in \cite{CKMS,FS,Hag,Kor,Mat2,Raf,Stu,Tak}.
For the numerical methods and comparison such as the finite difference time
domain (FDTD) methods \cite{BHM,Hamm,NSG},
time-splitting Fourier spectral (TSFP) methods \cite{BL,BZ,FSS,HJMSZ}
and Runge-Kutta discontinuous Galerkin methods \cite{WT,XST},
we refer to \cite{BL,BZ,BHM,FSS,Hamm,HL,HJMSZ,NSG}
and references therein.
However, for the NLDE (\ref{SDEd}) (or (\ref{SDEdd})) with $0<\varepsilon\ll 1$,
i.e. nonrelativistic limit regime (or the scaled speed of light goes to infinity),
the analysis and efficient computation of the NLDE (\ref{SDEd}) (or (\ref{SDEdd}))
are mathematically rather complicated issues.
The main difficulty is due to that the solution is highly oscillatory in time
and the corresponding energy functionals (\ref{engery60}) and (\ref{engery65})
are indefinite \cite{BMP, Esteban} and become unbounded when $\varepsilon\to0$.
For the Dirac equation, i.e. $\bF(\Psi)\equiv0$ in (\ref{nolc})
(or $\bF(\Phi)\equiv 0$ in (\ref{nolcc})),
there are extensive mathematical analysis of the (semi)-nonrelativistic
limits \cite{Hun,BMP,Gri,NM,White}.
For the NLDE (\ref{SDEd}) (or (\ref{SDEdd})), similar analysis of the nonrelativistic limits
has been done  in \cite{Foldy,N}.
These rigorous analytical results show that the solution
propagates waves with wavelength $O(\varepsilon^2)$ and
$O(1)$ in time and space, respectively, when $0<\varepsilon\ll 1$.
In fact, the oscillatory structure of the solution of the NLDE (\ref{SDEd}) (or (\ref{SDEdd}))
when $0<\varepsilon\ll 1$ can be formally observed from its dispersion relation
(\ref{disperg}) (or (\ref{dispergg})). To illustrate this further, Figure \ref{Oscillation}
shows the solution of the NLDE (\ref{SDEdd}) with $d=1$,
$V(t,x)=\frac{1-x}{1+x^2}$, $A_1(t,x)=\frac{(1+x)^2}{1+x^2}$, $\lambda_1=-1$, $\lambda_2=0$ and $\Phi_0(x)=\left(\exp(-x^2/2),\exp(-(x-1)^2/2)\right)^T$ for different $\varepsilon$.

The highly oscillatory nature of the solution of the NLDE (\ref{SDEd})  (or
(\ref{SDEdd})) causes severe numerical burdens in practical computation,
making the numerical approximation of the NLDE (\ref{SDEd})  (or
(\ref{SDEdd})) extremely challenging and costly in the
nonrelativistic regime $0<\varepsilon\ll 1$.
Recently, we compared the spatial/temporal resolution in term of
$\varepsilon$ and established rigorous error estimates of
the FDTD methods, TSFP methods for the Dirac equation in the nonrelativistic limit regime
\cite{BCJ}, and proposed a new uniformly accurate  multiscale time integrator pseudospectral method \cite{BCJ1}.
To our knowledge, so far there are few results on the
numerics of the NLDE in the nonrelativistic limit regime.
The aim of this paper is to study the efficiency of the Crank-Nicolson finite difference (CNFD)
and TSFP methods applied
to the NLDE in the nonrelativistic limit regime, to propose the
exponential wave integrator Fourier pseudospectral (EWI-FP) method and
to compare their resolution capacities in this regime. We
start with the detailed analysis on the  convergence of
the CNFD method by paying particular attention to
how the error bound depends explicitly on the small parameter $\eps$
in addition to the mesh size $h$ and time step $\tau$. Based on the
estimate, in order to obtain `correct' numerical approximations
when $0<\eps\ll1$, the meshing strategy requirement
($\eps$-scalability) for the CNFD  (and FDTD) is:
$h = O(\sqrt{\eps})$ and $\tau = O(\eps^3)$,
which suggests that the CNFD (and FDTD) is
computationally expensive for the NLDE (\ref{SDEd}) (or (\ref{SDEdd})) as $0<\eps\ll 1$. To
relax the $\eps$-scalability, we then propose the EWI-FP method and compare
it with the TSFP method, whose $\eps$-scalability are optimal for both time and space in view
of the inherent oscillatory nature. The key ideas of the EWI-FP
are: (i) to apply the Fourier pseudospectral  discretization for spatial derivatives;
and (ii) to adopt the exponential wave integrator (EWI), which was well demonstrated in the 
literatures that it has favorable properties
compared to standard time integrators for oscillatory  differential
equations \cite{Gautschi0,Gautschi-type-3}, for integrating
the ordinary differential equations (ODEs)
in phase space \cite{Gautschi0,Gautschi-type-3}.
Rigorous error
estimates show that the $\eps$-scalability of the EWI-FP method  is
$h = O(1)$, and $\tau= O(\eps^2)$
for the NLDE with  external
electromagnetic potentials,  meanwhile,
the $\eps$-scalability of the TSFP method
is $h = O(1)$ and $\tau= O(\eps^2)$.
Thus, the EWI-FP and TSFP offer compelling advantages
over CNFD (and FDTD) for the NLDE in temporal and spatial resolution when
$0<\eps\ll 1$. In particular, under suitable choices of the time step, the error estimates
of TSFP are much better than EWI-FP, which suggests that TSFP performs best in the nonrelativistic limit
regime.

The rest of this paper is organized as follows. In Section \ref{sec2},  the CNFD method
is reviewed and its convergence is analyzed in the nonrelativistic limit regime.
In Section \ref{sec3}, an EWI-FP
method is proposed and analyzed rigorously. In Section
\ref{sec4}, a TSFP  method
is reviewed and analyzed rigorously. In Section
\ref{sec5}, numerical comparison results are reported. Finally,
some concluding remarks are drawn in Section \ref{sec6}. The mathematical proofs of the error
estimates are given in the Appendices.
Throughout the paper, we adopt the standard notations of Sobolev spaces, use
the notation $p \lesssim q$ to represent that there exists a generic
constant $C>0$ which is independent of $h$, $\tau$ and $\eps$ such that
$|p|\le C\,q$.


\section{The Crank-Nicolson finite difference (CNFD) method}\label{sec2}
\setcounter{equation}{0}
\setcounter{table}{0}
\setcounter{figure}{0}

In this section, we apply the CNFD method to the
NLDE \eqref{SDEdd} (or (\ref{SDEd})) with external electromagnetic
field and analyze its conservation law and convergence in the
nonrelativistic limit regime.
For simplicity of notations, we shall only present the numerical method and its analysis
for (\ref{SDEdd}) in 1D. Generalization to (\ref{SDEd}) and/or higher dimensions is
straightforward and results remain valid without modifications.
Similar to most works in the
literatures for the analysis and computation of the NLDE
(cf. \cite{Alv1992, AKV, FS, HL, WT, XST} and references therein),
in practical computation, we truncate the whole space problem onto an
interval $\Omega=(a, b)$ with periodic boundary conditions, which is
large enough such that the truncation error is negligible.
In 1D, the NLDE \eqref{SDEdd} with periodic boundary conditions collapses to
\bea
\label{NLD1d}
&&i\partial_t\Phi=\left[-\frac{i}{\varepsilon}\sigma_1
\partial_x+\frac{1}{\varepsilon^2}\sigma_3\right]\Phi+\left[
V(t,x)I_2-A_1(t,x)\sigma_1+\bF(\Phi)\right]\Phi,
 \qquad x\in\Omega,\quad t>0,\\
&&\Phi(t,a)=\Phi(t,b),\quad\partial_x\Phi(t,a) =\partial_x \Phi(t,b),\quad t\geq 0, \qquad
\Phi(0,x) = \Phi_0(x),\quad x\in\overline{\Omega},
\label{NLD_I}
\eea
where $\Phi:=\Phi(t,x)\in {\mathbb C}^2$, $\Phi_0(a)=\Phi_0(b)$, $\Phi_0^{\prime}(a)=\Phi_0^{\prime}(b)$ and
$\bF(\Phi)$ is given in (\ref{nolcc}).

\subsection{The CNFD  method}

Choose mesh size $h:=\Delta x =\frac{b-a}{M}$ with $M$ being a  positive integer,
time step $\tau := \Delta t >0$ and denote the grid points and time steps as:
$$
x_j:=a+jh,\quad j=0, 1, \ldots ,M;\qquad t_n:= n\tau,\quad n=0, 1, 2,\ldots .
$$
Denote $X_M=\{\beU=(U_0,U_1,...,U_M)^T\ |\ U_j\in {\mathbb C}^2, j=0,1,\ldots,M, \ U_0=U_M\}$
and  we always use $U_{-1}=U_{M-1}$ if it is involved.
The standard $l^2$-norm in $X_M$ is given as
\begin{eqnarray}
\|\beU\|^2_{l^2}=h\sum^{M-1}_{j=0}|U_j|^2, \qquad  \beU\in X_M.
\end{eqnarray}
Let $\Phi_{j}^{n}$  be the numerical approximation of $\Phi(t_n,x_j)$ and
 $V_j^n=V(t_n,x_j)$, $V_j^{n+1/2}=V(t_n+\tau/2,x_j)$, $A_{1,j}^n=A_1(t_n,x_j)$,
 $A_{1,j}^{n+1/2}=A_1(t_n+\tau/2,x_j)$, $\bF_j^n=\bF(\Phi_j^n)$ and
$\bF_j^{n+1/2}=\frac{1}{2}\left[\bF(\Phi_j^n)+\bF(\Phi_j^{n+1})\right]$
 for $0\le j\le M$ and $n\ge0$.
Denote $\Phi^n=\left(\Phi_0^n, \Phi_1^n, \ldots, \Phi_M^n\right)^T\in X_M$ as
the solution vector at $t=t_n$.
Introduce the finite difference discretization operators for $j=0,1,\ldots, M-1$ and $n\ge0$ as:
\begin{equation*}
\delta_t^+\Phi_{j}^{n}=\frac{\Phi_{j}^{n+1}-\Phi_{j}^{n}}{\tau},\qquad
\delta_t\Phi_{j}^{n}=\frac{\Phi_{j}^{n+1}-\Phi_{j}^{n-1}}{2\tau},\qquad\delta_x\Phi_{j}^{n}
=\frac{\Phi_{j+1}^{n}-\Phi_{j-1}^{n}}{2h},\qquad
\Phi_{j}^{n+\frac{1}{2}}=\frac{\Phi_{j}^{n+1}+\Phi_{j}^{n}}{2}.
\end{equation*}

Here we consider the CNFD method
 to discretize the NLDE \eqref{NLD1d}:

\begin{equation}
\label{nlcnfd1}
i\delta_t^+\Phi_{j}^{n}=\left[-\frac{i}{\vep}\sigma_1\delta_x +\frac{1}{\vep^2}\sigma_3
+V_j^{n+1/2}I_2-A_{1,j}^{n+1/2}\sigma_1+\bF^{n+1/2}_j\right]\Phi_{j}^{n+1/2},
\quad 0\le j\le M-1, \quad n\ge0.
\end{equation}

\bigskip

\noindent The initial and boundary conditions in (\ref{NLD_I}) are discretized as:
\begin{eqnarray}
\label{nlbdc}
\Phi_{M}^{n+1}=\Phi_{0}^{n+1},\quad \Phi_{-1}^{n+1}=\Phi_{M-1}^{n+1},\quad n\geq 0,
\qquad \Phi_{j}^{0}=\Phi_0(x_j),\quad j=0,1,...,M.
\end{eqnarray}

\subsection{Conservation law and error estimates}
Let  $0<T<T^*$ with $T^*$ being the maximal existence time of
the solution, and denote $\Omega_{T}=[0,T]\times \Omega$.
We assume the electromagnetic potentials
$V\in C(\overline{\Omega}_T)$ and $A_1\in C(\overline{\Omega}_T)$ and denote
\[
(A) \hskip4cm   V_{\rm max}:=\max_{(t,x)\in\overline{\Omega}_T}|V(t,x)|,
\qquad A_{1,\rm max}:=\max_{(t,x)\in\overline{\Omega}_T}|A_1(t,x)|. \hskip5cm
\]

For the CNFD method \eqref{nlcnfd1}, similar to
the mass and energy conservation of the Dirac equation  \cite{BCJ}, 
we have the following conservative properties, of which the proof is omitted here for brevity.

\begin{lemma}
\label{Mass_STA_NLD_FDTD}
The  CNFD method \eqref{nlcnfd1} conserves
the mass in the discretized level, i.e.
\be\label{nlcons976}
\|\Phi^n\|_{l^2}^2:=h\sum_{j=0}^{M-1}|\Phi_j^{n}|^2\equiv  h\sum_{j=0}^{M-1}|\Phi_j^{0}|^2
=\|\Phi^0\|_{l^2}^2=h\sum_{j=0}^{M-1}|\Phi_0(x_j)|^2,\qquad n\ge0.
\ee
Furthermore, if $V(t,x)=V(x)$ and $A_1(t,x)=A_1(x)$ are time independent, the
CNFD method \eqref{nlcnfd1} conserves the energy as well
\bea\label{nlcons:energy}
E_h^n&=&h\sum\limits_{j=0}^{M-1}\left[-\frac{i}{\vep}(\Phi_j^{n})^*\sigma_1\delta_x\Phi_j^{n}
+\frac{1}{\vep^2}(\Phi_j^{n})^*\sigma_3\Phi_j^{n}+V(x_j)|\Phi_j^{n}|^2-
A_{1}(x_j)(\Phi_j^{n})^*\sigma_1\Phi_j^{n}+G(\Phi_j^n)\right]\nonumber\\
&\equiv&E_h^0,\qquad n\ge0,
\eea
where $G(\Phi)$ is given in (\ref{GPhi}). 
\end{lemma}

Next, we consider the error analysis of the CNFD \eqref{nlcnfd1}. Motivated by the analytical results
of the NLDE, we assume that the exact solution of \eqref{NLD1d} satisfies $\Phi\in C^{3}([0,T]; (L^{\infty}(\Omega))^2)\cap C^{2}([0,T]; (W_p^{1,\infty}(\Omega))^2)
\cap C^{1}([0,T]; (W_p^{2,\infty}(\Omega))^2)\cap C([0,T]; (W_p^{3,\infty}(\Omega))^2)$ and
\[
(B) \hskip3cm \left\|\frac{\partial^{r+s}}{\partial t^{r}
\partial x^{s}}\Phi\right\|_{L^\infty([0,T];(L^{\infty}(\Omega))^2)}
\lesssim \frac{1}{\varepsilon^{2r}},
\quad0\leq r\leq 3,\ 0\leq r+s\leq 3, \quad 0<\vep\le 1,\hskip4cm
\]
where $W_p^{m,\infty}(\Omega)=\{u\ |\ u\in W^{m,\infty}(\Omega),\  \partial_x^l u(a)=\partial_x^l u(b),\
l=0,\ldots,m-1\}$ for $m\ge1$ and here the boundary values are understood in the trace sense. In the subsequent discussion,
we will omit $\Omega$ when referring to the space norm taken on $\Omega$. We denote
\be
M_0:=\max_{0\le t\le T}\|\Phi(t,x)\|_{L^{\infty}}\lesssim 1.
\ee

Define the grid error function $\bee^n=(\bee_0^n,\bee_1^n,\ldots,\bee_M^n)^T\in X_M$ as:
\begin{equation}
\label{nlerr}
\bee_{j}^n = \Phi(t_n,x_j) - \Phi_{j}^n, \qquad j=0,1,\ldots,M, \quad n\ge0,
\end{equation}
with $\Phi_j^n$ being the approximations obtained from the CNFD method.

For the CNFD \eqref{nlcnfd1}, we can establish the error bound (see its proof in Appendix A).

\begin{theorem}
\label{nlthm_cnfd}
Assume $0<\tau\lesssim \vep^3h^{\frac{1}{4}}$,
under the assumptions $(A)$ and $(B)$, there exist constants $h_{0}>0$ and $\tau_0>0$
sufficiently small and independent of $\vep$, such that for any
$0<\vep\leq1$, when $0<h\leq h_{0}$ and $0<\tau\leq\tau_0$ satisfying $0<h\lesssim \vep^{\frac{2}{3}}$,
we have the following error estimate for the CNFD method \eqref{nlcnfd1} with \eqref{nlbdc}
\begin{equation} \label{nlebcnfd97}
\|\bee^{n}\|_{l^{2}}\lesssim \frac{h^{2}}{\varepsilon}+\frac{\tau^{2}}{\varepsilon^{6}},\qquad \|\Phi^n\|_{l^{\infty}}\leq 1+M_0,
\qquad 0\leq n\leq\frac{T}{\tau}.
\end{equation}
\end{theorem}

\begin{remark}
The above Theorem is still valid in higher dimensions provided that the conditions
$0<\tau\lesssim\vep^3h^{\frac{1}{4}}$ and
$0<h\lesssim \vep^{\frac{2}{3}}$ are replaced by $0<\tau\lesssim \vep^3h^{C_d}$ and
$0<h\lesssim \vep^{\frac{1}{2(1-C_d)}}$, respectively, with $C_d=\frac{d}{4}$ for $d=1,2,3$.
\end{remark}

\begin{remark}\label{remark:1}
Similar to the Dirac equation, we can easily extend other
 finite difference time domain (FDTD) methods including the leap-frop finite difference (LFFD) and
 semi-implicit finite difference (SIFD) \cite{BCJ,Jia} to the NLDE, and the error bounds are the same as
 those in Theorem \ref{nlthm_cnfd}.
\end{remark}

Based on Theorem \ref{nlthm_cnfd} and Remark \ref{remark:1}, the  CNFD method (and FDTD methods)
has the following temporal/spatial resolution capacity for the NLDE in the nonrelativistic limit regime. In fact,
given an accuracy bound $\delta_0>0$, the $\vep$-scalability of the CNFD method  is:
\begin{equation*}
\tau=O\left(\vep^3 \sqrt{\delta_0}\right)=O(\vep^3), \qquad
h=O\left(\sqrt{\delta_0\vep}\right)=O\left(\sqrt{\vep}\right),
\qquad 0<\vep\ll1.
\end{equation*}


\section{An EWI-FP method and its analysis} \label{sec3}
\setcounter{equation}{0}
\setcounter{table}{0}
\setcounter{figure}{0}

In this section, we propose an EWI-FP
method  to solve the NLDE (\ref{NLD1d})
and establish its error bound.

\subsection{The EWI-FP method}

Denote $\mu_l=\frac{2l\pi}{b-a}$ for $l\in{\mathbb Z}$ and
\[
Y_{M}=Z_M\times Z_M, \quad
Z_M={\rm span}\left\{\phi_l(x)=e^{i\mu_l(x-a)},\ l=-\frac{M}{2}, -\frac{M}{2}+1,\ldots, \frac{M}{2}-1\right\}.
\]
Let $[C_p(\overline{\Omega})]^2$ be the function space consisting of all periodic vector function $U(x):\ \overline{\Omega}=[a,b]\to {\mathbb C}^2$. For any $U(x)\in [C_p(\overline{\Omega})]^2$ and $U\in X_M$,
define $P_M:\ [L^2(\Omega)]^2\rightarrow Y_M$ as the standard projection operator \cite{ST}, $I_M:\ [C_p(\overline{\Omega})]^2\rightarrow Y_M$ and $I_M:X_M\rightarrow Y_M$ as the standard interpolation operator, i.e.
\begin{equation*}
(P_MU)(x)=\sum_{l=-M/2}^{M/2-1}\widehat{U}_l\, e^{i\mu_l(x-a)},\quad (I_MU)(x)
=\sum_{l=-M/2}^{M/2-1}\widetilde{U}_l\, e^{i\mu_l(x-a)},\qquad
a\leq x\leq b,
\end{equation*}
with
\be\label{nlfouriercoef}
\widehat{U}_l=\frac{1}{b-a}\int_a^bU(x)\,e^{-i\mu_l(x-a)}\,dx,\quad
\widetilde{U}_l=\frac{1}{M}\sum_{j=0}^{M-1}U_j\, e^{-2ijl\pi/M},\quad
l=-\frac{M}{2}, \ldots, \frac{M}{2}-1,
\ee
where $U_j=U(x_j)$ when $U$ is a function.

The Fourier spectral discretization for the NLDE (\ref{NLD1d}) is as follows:

\noindent Find $\Phi_M:=\Phi_{M}(t,x)\in Y_M$, i.e.
\begin{eqnarray}
\label{nlFPR}
\Phi_{M}(t,x)=\sum^{M/2-1}_{l=-M/2}\widehat{(\Phi_M)}_l(t)\,
e^{i\mu_l(x-a)},\qquad a\le x\le b, \qquad t\ge0,
\end{eqnarray}
such that
\bea
\label{nlFTD}
i\partial_t\Phi_{M}=\left[-\frac{i}{\varepsilon}\sigma_1\partial_x
+\frac{1}{\varepsilon^2}\sigma_3\right]\Phi_{M}+P_M\left[\left(
V(t,x)I_2-A_1(t,x)\sigma_1+\bF(\Phi_M)\right)\Phi_M\right], \ \ a<x<b,\ \ t>0.\
\eea
Substituting (\ref{nlFPR}) into (\ref{nlFTD}), noticing the orthogonality of $\phi_l(x)$,
we get
\begin{equation}
i\,\frac{d}{dt}\widehat{(\Phi_M)}_{l}(t)=\left[\frac{\mu_l}{\varepsilon}\sigma_1
+\frac{1}{\varepsilon^2}\sigma_3\right]\widehat{(\Phi_M)}_{l}(t)+
\widehat{\bG(\Phi_M)}_l(t),\quad t>0,\quad l=-\frac{M}{2}, -\frac{M}{2}+1, \ldots, \frac{M}{2}-1,
\end{equation}
where
\be\label{bGtx86}
\bG(\Phi_M)=\left(V(t,x)I_2-A_1(t,x)\sigma_1+\bF(\Phi_M)\right)\Phi_M,\qquad x\in\Omega, \quad t\ge0.
\ee
For $t\ge0$ and when it is near $t=t_n$ $(n\geq 0)$, we rewrite the above ODEs  as
\be\label{nlODE765}
i\,\frac{d}{ds}\widehat{(\Phi_M)}_{l}(t_n+s)=\frac{1}{\vep^2}\Gamma_l\,
\widehat{(\Phi_M)}_{l}(t_n+s)+\widehat{\bG(\Phi_M)_l^n}(s),\quad s>0,\quad
l=-\frac{M}{2}, -\frac{M}{2}+1, \ldots, \frac{M}{2}-1,
\ee
where $\Gamma_l=\mu_l\vep\sigma_1+\sigma_3=Q_l\, D_l\, (Q_l)^*$ with $\delta_l=\sqrt{1+\vep^2\mu_l^2}$
and
\begin{equation}\label{nleq:Gamma}
\Gamma_l=\begin{pmatrix}
1 &\mu_l\vep\\
\mu_l\vep &-1
\end{pmatrix}, \quad Q_l=\begin{pmatrix}
\frac{1+\delta_l}{\sqrt{2\delta_l(1+\delta_l)}} &-\frac{\vep\mu_l}{\sqrt{2\delta_l(1+\delta_l)}}\\
\frac{\vep\mu_l}{\sqrt{2\delta_l(1+\delta_l)}} &\frac{1+\delta_l}{\sqrt{2\delta_l(1+\delta_l)}}
\end{pmatrix}, \quad D_l=\begin{pmatrix}
\delta_l &0\\
0 &-\delta_l\\
\end{pmatrix},
\end{equation}
and
\begin{equation}\label{nleq:fdef}
\widehat{\bG(\Phi_M)_l^n}(s)=\widehat{\bG(\Phi_M)}_l(t_n+s),\qquad s\ge0, \quad n\ge0.
\end{equation}

Solving the above ODE (\ref{nlODE765}) via the integration factor method, we obtain
\begin{eqnarray}
\label{nlEWIS}
\widehat{(\Phi_M)}_l(t_n+s)=e^{-is\Gamma_l/\vep^2}\widehat{(\Phi_M)}_l(t_n)
-i\int_0^{s}e^{i(w-s)\Gamma_l/\vep^2}\widehat{\bG(\Phi_M)_l^n}(w)\,dw, \qquad s\ge0.
\end{eqnarray}
Taking $s=\tau$ in (\ref{nlEWIS}) we have
\begin{equation}\label{nlEWIn}
\widehat{(\Phi_M)}_l(t_{n+1})=e^{-i\tau\Gamma_l/\vep^2}\widehat{(\Phi_M)}_l(t_{n})
-i\int_0^{\tau}e^{\frac{i(w-\tau)}{\vep^2}\Gamma_l}
\widehat{\bG(\Phi_M)_l^n}(w)dw.
\end{equation}
To obtain a numerical method with second  order accuracy in time, we approximate the integral in \eqref{nlEWIn}
via the Gautschi-type rule which has been widely
used  for integrating highly oscillatory ODEs \cite{BD, Gautschi0, Gautschi-type-3}:
\bea\label{nleq:inteapp0}
\int_0^{\tau}e^{\frac{i(w-\tau)}{\vep^2}\Gamma_l}\widehat{\bG(\Phi_M)_l^0}(w)\,dw
\approx\int_0^{\tau}e^{\frac{i(w-\tau)}{\vep^2}\Gamma_l}\,dw\,\widehat{\bG(\Phi_M)_l^0}(0)
=-i\vep^2\Gamma_l^{-1}\left[I_2-e^{-\frac{i\tau}{\vep^2}\Gamma_l}\right]
\widehat{\bG(\Phi_M)_l^0}(0), \quad
\eea
and for $n\ge 1$
\begin{align}\label{nleq:inteappn}
&\int_0^{\tau}e^{\frac{i(w-\tau)}{\vep^2}\Gamma_l}\widehat{\bG(\Phi_M)_l^n}(w)
\approx\int_0^{\tau}e^{\frac{i(w-\tau)}{\vep^2}\Gamma_l}
\left(\widehat{\bG(\Phi_M)_l^n}(0)+w\,\delta_t^-\widehat{\bG(\Phi_M)_l^n}(0)\right)dw\nonumber\\
&\quad =-i\vep^2\Gamma_l^{-1}\left[I_2-e^{-\frac{i\tau}{\vep^2}\Gamma_l}\right]
\widehat{\bG(\Phi_M)_l^n}(0)+\left[-i\vep^2\tau\Gamma_l^{-1}+\vep^4\Gamma_l^{-2}
\left(I_2-e^{-\frac{i\tau}{\vep^2}\Gamma_l}\right)\right]\delta_t^-\widehat{\bG(\Phi_M)_l^n}(0),
\end{align}
where we have approximated the time derivative $\partial_t\widehat{\bG(\Phi_M)_l^n}(s)$
at $s=0$ by finite difference as
\be
\partial_t\widehat{\bG(\Phi_M)_l^n}(0)\approx\delta_t^-\widehat{\bG(\Phi_M)_l^n}(0)
=\frac{1}{\tau}\left[\widehat{\bG(\Phi_M)_l^{n}}(0)-\widehat{\bG(\Phi_M)_l^{n-1}}(0)\right].
\ee
Now, we are ready to describe our scheme. Let $\Phi_M^n(x)$ be the approximation of  $\Phi_{M}(t_n,x)$ ($n\ge0$).
Choosing $\Phi_M^0(x)=(P_M\Phi_0)(x)$,  an {\sl exponential wave integrator  Fourier spectral} (EWI-FS)
discretization for the NLDE \eqref{NLD1d} is to update the numerical approximation $\Phi^{n+1}_M(x)\in Y_M$ ($n=0,1,\ldots$) as
\begin{equation}\label{nlAL1}
\Phi_M^{n+1}(x)=
\sum_{l=-M/2}^{M/2-1}\widehat{(\Phi_M^{n+1})}_l\, e^{i\mu_l(x-a)},\qquad a\leq x\leq b,\qquad n \ge0,
\end{equation}
where for $l = -\frac{M}{2}, .., \frac{M}{2}-1$,
\be\label{nlCoe2}
\widehat{(\Phi_M^{n+1})}_l=\left\{\ba{ll}
e^{-i\tau\Gamma_l/\vep^2}\widehat{(\Phi_M^0)}_l
-\vep^2\Gamma_l^{-1}\left[I_2-e^{-\frac{i\tau}{\vep^2}\Gamma_l}\right]
\widehat{\bG(\Phi^0_M)}_l, &n=0,\\
e^{-i\tau\Gamma_l/\vep^2}\widehat{(\Phi_M^n)}_l
-iQ_l^{(1)}(\tau)
\, \widehat{\bG(\Phi^n_M)}_l-iQ_l^{(2)}(\tau)\delta_t^-\widehat{\bG(\Phi^n_M)}_l, &n\ge1,\\
\ea\right.
\ee
with the matrices $Q_l^{(1)}(\tau)$ and $Q_l^{(2)}(\tau)$ given as
\be
Q_l^{(1)}(\tau)=-i\vep^2\Gamma_l^{-1}\left[I_2-e^{-\frac{i\tau}{\vep^2}\Gamma_l}\right],
\quad
Q_l^{(2)}(\tau)=-i\vep^2\tau\Gamma_l^{-1}+\vep^4\Gamma_l^{-2}\left(I_2-e^{-\frac{i\tau}{\vep^2}\Gamma_l}\right),
\ee
and
\be\label{bGphi789}
\bG(\Phi^n_M)=\left(V(t_n,x)I_2-A_1(t_n,x)\sigma_1+\bF(\Phi_M^n)\right)\Phi_M^n,\qquad n\ge0.
\ee

The above procedure is not suitable in practice due to the difficulty in computing the Fourier coefficients through integrals in \eqref{nlfouriercoef}.
 Here we present an efficient implementation by choosing $\Phi_M^0(x)$ as the interpolant of $\Phi_0(x)$ on the grids $\left\{x_j, j = 0, 1,\ldots, M\right\}$ and approximate the integrals
in (\ref{nlfouriercoef}) by a quadrature rule.

Let $\Phi_j^n$  be the numerical approximation of $\Phi(t_n,x_j)$ for $j = 0, 1, 2, \ldots, M$ and $n\geq 0$, and denote $\Phi^n\in X_M$ as the vector with components $\Phi_j^n$.
 Choosing $\Phi^0_j=\Phi_0(x_j)$ ($j = 0, 1,\ldots, M$), an {\sl EWI Fourier pseudospectral} (EWI-FP) method for computing $\Phi^{n+1}$ for $n\geq 0$ reads
\begin{equation}\label{nlAL4}
\Phi_{j}^{n+1} = \sum_{l=-M/2}^{M/2-1}\widetilde{(\Phi^{n+1})}_le^{2ijl\pi/M}, \quad j = 0, 1, ..., M,
\end{equation}
where
\begin{equation}\label{nlAL5}
\widetilde{(\Phi^{n+1})}_l=\left\{\ba{ll}
e^{-i\tau\Gamma_l/\vep^2}\widetilde{(\Phi^0)}_l-\vep^2\Gamma_l^{-1}
\left[I_2-e^{-\frac{i\tau}{\vep^2}\Gamma_l}\right]\widetilde{\bG(\Phi^0)}_l, &n=0,\\
e^{-i\tau\Gamma_l/\vep^2}\widetilde{(\Phi^{n})}_l-iQ_l^{(1)}(\tau)
\,\widetilde{\bG(\Phi^{n})}_l-iQ_l^{(2)}(\tau)\delta_t^-\widetilde{\bG(\Phi^n)}_l, &n\ge1.\\
\ea\right.
\end{equation}
The EWI-FP (\ref{nlAL4})-(\ref{nlAL5}) is explicit, and can be solved efficiently by the
fast Fourier transform (FFT). The memory
cost is $O(M)$ and the computational cost per time step is $O(M \ln M)$.

  Similar to the analysis of the EWI-FP method for the Dirac equation in \cite{BCJ},
we can obtain that the EWI-FP for the NLDE is stable under the stability condition (details
are omitted here for brevity)
\begin{equation} \label{nlsbcexp}
0<\tau\lesssim 1, \quad 0<\vep\le 1.
\end{equation}

\subsection{An error estimate}

In order to obtain an error estimate for the EWI methods (\ref{nlAL1})-(\ref{nlCoe2}) and (\ref{nlAL4})-(\ref{nlAL5}),
motivated by the results in \cite{Foldy,N}, we assume that there exists an integer $m_0\geq 2$ such that the exact solution $\Phi(t,x)$ of the NLDE (\ref{NLD1d}) satisfies
\begin{align*}
(C)\hskip2cm \|\Phi\|_{ L^\infty([0,T]; (H_p^{m_0})^2)}\lesssim 1,\quad
\|\partial_t\Phi\|_{L^{\infty}([0, T];(L^2)^2)}\lesssim \frac{1}{\varepsilon^2},\quad \|\partial_{tt}\Phi\|_{L^{\infty}([0,T];(L^2)^2)}\lesssim\frac{1}{\varepsilon^4},\hskip4cm
\end{align*}
where $H^k_p(\Omega)=\{u\ |\ u\in H^{k}(\Omega),\  \partial_x^l u(a)=\partial_x^l u(b),\
l=0,\ldots,k-1\}$.
In addition, we assume electromagnetic potentials satisfy
\begin{equation*}
(D)\hskip5cm \|V\|_{W^{2,\infty}([0,T];L^\infty)}+\|A_1\|_{W^{2,\infty}([0,T];L^\infty)}\lesssim 1.\hskip7cm
\end{equation*}

We can establish the following error estimate for the EWI-FS method (see its proof in Appendix B).
\begin{theorem}
\label{nlthm_EWI}
Let $\Phi_M^n(x)$ be the approximation obtained from the EWI-FS (\ref{nlAL1})-(\ref{nlCoe2}). Assume $0<\tau\lesssim\varepsilon^2h^{1/4}$, under the
assumptions $(C)$ and $(D)$, there exist $h_0>0$ and $\tau_0>0$ sufficiently small
and independent of $\varepsilon$ such that, for any $0<\varepsilon\leq1$,
when $0<h\leq h_0$ and $0<\tau\leq\tau_0$, we have the error estimate
\begin{equation}
\label{nlthm_eq_EWI}
\|\Phi(t_n,x)-\Phi_M^n(x)\|_{L^2}\lesssim\frac{\tau^2}{\varepsilon^4}+h^{m_0},
\qquad \|\Phi_M^n(x)\|_{L^{\infty}}\leq 1+M_0,\qquad 0\leq n\leq\frac{T}{\tau}.
\end{equation}
\end{theorem}
\begin{remark}
The same error estimate in Theorem \ref{nlthm_EWI} holds for the
EWI-FP (\ref{nlAL4})-(\ref{nlAL5}) and the proof is quite similar to that of Theorem \ref{nlthm_EWI}.
In addition, the above Theorem is still valid in higher dimensions provided that the condition
$0<\tau\lesssim \vep^2h^{\frac{1}{4}}$ is replaced by $0<\tau\lesssim \vep^2 h^{C_d}$.
\end{remark}

From this theorem, the temporal/spatial resolution capacity of the EWI-FP method for the NLDE
in the nonrelativistic limit regime is: $h=O(1)$ and $\tau=O(\vep^2)$. In fact,
for a given accuracy bound $\delta_0>0$, the $\vep$-scalability of the EWI-FP is:
\begin{equation*}
\tau=O\left(\vep^2 \sqrt{\delta_0}\right)=O(\vep^2), \qquad  h=O\left(\delta_0^{1/m_0}\right)=O\left(1\right),
\qquad 0<\vep\ll1.
\end{equation*}

Similar to the Appendix D in \cite{BCJ} for the Dirac equation,
it is straightforward to generalize the EWI-FP to the NLDE (\ref{SDEdd}) in 2D and
(\ref{SDEd}) in 1D, 2D and 3D and the details are omitted here for brevity.


\section{A TSFP method and its optimal error bounds} \label{sec4}
\setcounter{equation}{0}
\setcounter{table}{0}
\setcounter{figure}{0}
In this section, we present a time-splitting Fourier pseudospectral (TSFP)
method for the NLDE \eqref{NLD1d}.

\subsection{The TSFP method}
From time $t = t_n$ to time $t = t_{n+1}$,
the NLDE \eqref{NLD1d} is split into two steps. One solves first
\be
\label{nldir1st1}
i\partial_t\Phi(t,x)=\left[-\frac{i}{\varepsilon}\sigma_1
\partial_x+\frac{1}{\varepsilon^2}\sigma_3\right]\Phi(t,x), \quad x\in\Omega,
\ee
with the periodic boundary condition \eqref{NLD_I}
for the time step of length $\tau$, followed by solving
\be
\label{nldir1st2}
i\partial_t\Phi(t,x)=\left[
V(t,x)I_2-A_1(t,x)\sigma_1+\bF(\Phi(t,x))\right]\Phi(t,x), \quad x\in\Omega,
\ee
for the same time step.
Eq. \eqref{nldir1st1} will be first discretized
in space by the Fourier spectral method and then
integrated (in phase or Fourier space) in
time {\sl exactly} \cite{BCJ,BL}. For the ODEs \eqref{nldir1st2},
multiplying  $\Phi^*(t,x)$ from the left, we get
\be\label{Phisp12}
i\Phi^*(t,x)\partial_t\Phi(t,x)=\Phi^*(t,x)\left[
V(t,x)I_2-A_1(t,x)\sigma_1+\bF(\Phi(t,x))\right]\Phi(t,x), \quad x\in\Omega.
\ee
Taking conjugate to both sides of the above equation, noticing (\ref{nolcc}), we obtain
\be \label{Phisp13}
-i\partial_t\Phi^*(t,x)\,\Phi(t,x)=\Phi^*(t,x)\left[
V(t,x)I_2-A_1(t,x)\sigma_1^*+\bF(\Phi(t,x))\right]\Phi(t,x), \quad x\in\Omega,
\ee
where $\sigma_1^*=\overline{\sigma_1}^T$.
Subtracting (\ref{Phisp13}) from (\ref{Phisp12}), noticing (\ref{nolcc}),
 $\sigma_1^*=\sigma_1$ and $\sigma_3^*=\sigma_3$ , we obtain for
$\rho(t,x)=|\Phi(t,x)|^2$
\be
\partial_t \rho(t,x)=0, \qquad t_n\le t\le t_{n+1}, \qquad x\in \Omega,
\ee
which immediately implies $\rho(t,x)=\rho(t_n,x)$.

If $A_1(t,x)\equiv 0$,
multiplying  \eqref{nldir1st2}
from left by $\Phi^*(t,x)\sigma_3$ and by a similar procedure, we get
$\Phi^*(t,x)\sigma_3\Phi(t,x)=  \Phi^*(t_n,x)\sigma_3\Phi(t_n,x)$ for
$t_n\le t\le t_{n+1}$ and $x\in\Omega$. Thus if $\lambda_1=0$ or $A_1(t,x)\equiv 0$, we have
\be\label{bfphitx}
\bF(\Phi(t,x))=\bF(\Phi(t_n,x)), \qquad t_n\le t\le t_{n+1}, \qquad x\in \Omega.
\ee
Plugging (\ref{bfphitx}) into (\ref{nldir1st2}), we obtain
\be
\label{nldir1st6}
i\partial_t\Phi(t,x)=\left[
V(t,x)I_2-A_1(t,x)\sigma_1+\bF(\Phi(t_n,x))\right]\Phi(t,x), \quad x\in\Omega,
\ee
which can be integrated {\sl analytically} in time as
\be
\Phi(t,x)=e^{-i\int_{t_n}^t \left[V(s,x)I_2-A_1(s,x)\sigma_1+\bF(\Phi(t_n,x))\right]ds}\,\Phi(t_n,x),
\qquad a\le x\le b, \quad t_n\le t\le t_{n+1}.
\ee
 In practical computation, if $\lambda_1=0$ or $A_1(t,x)\equiv 0$, from time $t=t_n$ to $t=t_{n+1}$, we
often combine the splitting steps via the Strang splitting \cite{St}
-- which results in a second order TSFP
method as
\begin{equation}\label{nleq:tsfp}
\begin{split}
&\Phi_j^{(1)}=\sum_{l=-M/2}^{M/2-1} e^{-i\tau \Gamma_l/2\vep^2}\,\widetilde{(\Phi^n)}_l\, e^{i\mu_l(x_j-a)}
=\sum_{l=-M/2}^{M/2-1} Q_l\, e^{-i \tau D_l/2\vep^2}\,(Q_l)^*\,
\widetilde{(\Phi^n)}_l\, e^{\frac{2ijl\pi}{M}},\\
&\Phi_j^{(2)}=e^{-i\int_{t_n}^{t_{n+1}} \bW(t,x_j)\,dt}\,\Phi_j^{(1)} = P_j\, e^{-i\, \Lambda_j}\,
P_j^*\, \Phi_j^{(1)},\qquad j=0,1,\ldots,M,\qquad n\ge0,\\
&\Phi_j^{n+1}=\sum_{l=-M/2}^{M/2-1} e^{-i \tau \Gamma_l/2\vep^2}\,\widetilde{(\Phi^{(2)})}_l\, e^{i\mu_l(x_j-a)}
=\sum_{l=-M/2}^{M/2-1} Q_l\, e^{-i \tau D_l/2\vep^2}\,(Q_l)^*\,
\widetilde{(\Phi^{(2)})}_l\, e^{\frac{2ijl\pi}{M}},
\end{split}
\end{equation}
where
$\int_{t_n}^{t_{n+1}}\bW(t,x_j)dt=V_j^{(1)}\, I_2-
A_{1,j}^{(1)}\,\sigma_1+\tau \bF(\Phi_j^{(1)})=(V_j^{(1)}+\tau \lambda_2|\Phi_j^{(1)}|^2) I_2-
A_{1,j}^{(1)}\,\sigma_1+\tau \lambda_1\, (\Phi_j^{(1)})^*\sigma_3\Phi_j^{(1)}\sigma_3$ $
 =P_j\, \Lambda_j\, P_j^*$ with $V_j^{(1)}=\int_{t_n}^{t_{n+1}}V(t,x_j)dt$,
$A_{1,j}^{(1)}=\int_{t_n}^{t_{n+1}}A_1(t,x_j)dt$, $\Lambda_j={\rm diag}(\Lambda_{j,+},\Lambda_{j,-})$, and
$\Lambda_{j,\pm}=V_j^{(1)}+\tau \lambda_2|\Phi_j^{(1)}|^2\pm \tau\lambda_1
(\Phi_j^{(1)})^*\sigma_3\Phi_j^{(1)}$ and $P_j=I_2$ if $A_{1,j}^{(1)}=0$, and resp.,
$\Lambda_{j,\pm}=V_j^{(1)}+\tau \lambda_2|\Phi_j^{(1)}|^2\pm A_{1,j}^{(1)}$ and
\be\label{eq:Gjn}
 P_j=P^{(0)}:=\begin{pmatrix}
\frac{1}{\sqrt{2}} &\frac{1}{\sqrt{2}}\\
-\frac{1}{\sqrt{2}} &\frac{1}{\sqrt{2}}\end{pmatrix},
\ee
if $A_{1,j}^{(1)}\ne 0$ and $\lambda_1=0$.

  Of course, if $\lambda_1\ne0$ and $A_1(t,x)\ne0$, then $\Phi^*(t,x)\sigma_3\Phi(t,x)$ is no longer
time-independent in the second step (\ref{nldir1st2}) due to the fact that $\sigma_1^*\sigma_3^*=\sigma_1\sigma_3\ne
\sigma_3\sigma_1$. In this situation, we will spit (\ref{nldir1st2}) into two steps as:
one first solves
\be
\label{nldir1st3}
i\partial_t\Phi(t,x)=\left[
V(t,x)I_2-A_1(t,x)\sigma_1\right]\Phi(t,x), \quad x\in\Omega,
\ee
followed by solving
\be
\label{nldir1st4}
i\partial_t\Phi(t,x)=\bF(\Phi(t,x))\,\Phi(t,x), \quad x\in\Omega.
\ee
Similar to the Dirac equation \cite{BCJ}, Eq. (\ref{nldir1st3}) can be
integrated {\sl analytically} in time. For Eq. (\ref{nldir1st4}),
both $\rho(t,x)$ and  $\Phi^*(t,x)\sigma_3\Phi(t,x)$ are invariant in time, i.e.
$\rho(t,x)\equiv \rho(t_n,x)$ and $\Phi^*(t,x)\sigma_3\Phi(t,x)\equiv
\Phi^*(t_n,x)\sigma_3\Phi(t_n,x)$ for $t_n\le t\le t_{n+1}$ and $x\in\bar\Omega$.
Thus it collapses to
\be
\label{nldir1st5}
i\partial_t\Phi(t,x)=\bF(\Phi(t_n,x))\,\Phi(t,x), \quad x\in\Omega,
\ee
and it can be integrated {\sl analytically} in time too. Similarly,
a second-order TSFP method can be designed provided that we replace
$\Phi^{(2)}$ in the third step by $\Phi^{(4)}$ and the second
step in (\ref{nleq:tsfp}) by
\begin{equation}\label{nleq:tsfp98}
\begin{split}
&\Phi_j^{(2)}=e^{-\frac{i}{2}\int_{t_n}^{t_{n+1}}
\bF(\Phi(t_n,x_j))\,dt}\,\Phi_j^{(1)} = e^{-i \,\Lambda_j^{(1)}}\, \Phi_j^{(1)},\\
&\Phi_j^{(3)}=e^{-i\int_{t_n}^{t_{n+1}}\left[
V(t,x_j)I_2-A_1(t,x_j)\sigma_1\right] \,dt}\,\Phi_j^{(2)} = P_j\, e^{-i\, \Lambda_j^{(2)}}\,
P_j^*\, \Phi_j^{(2)},\\
&\Phi_j^{(4)}=e^{-\frac{i}{2}\int_{t_n}^{t_{n+1}} \bF(\Phi(t_n,x_j))\,dt}\,\Phi_j^{(3)}
=  e^{-i \,\Lambda_j^{(1)}}\, \Phi_j^{(3)},\qquad j=0,1,\ldots,M,\qquad n\ge0,
\end{split}
\end{equation}
where $\Lambda_j^{(1)}={\rm diag}(\Lambda_{j,+}^{(1)},\Lambda_{j,-}^{(1)})$ with
$\Lambda_{j,\pm}^{(1)}=\frac{\tau}{2}\left[\lambda_2|\Phi_j^{(1)}|^2\pm \lambda_1
(\Phi_j^{(1)})^*\sigma_3\Phi_j^{(1)}\right]$,
$\Lambda_j^{(2)}={\rm diag}(\Lambda_{j,+}^{(2)},\Lambda_{j,-}^{(2)})$ with
$\Lambda_{j,\pm}^{(2)}=V_j^{(1)}\pm A_{1,j}^{(1)}$, and $P_j=I_2$ if $A_{1,j}^{(1)}=0$, and resp.,
$P_j=P^{(0)}$ if $A_{1,j}^{(1)}\ne 0$ for $j=0,1,\ldots,M$.

\begin{remark}
If the above definite integrals  cannot be evaluated analytically, we can
evaluate them numerically via the Simpson's quadrature rule as
\begin{align*}&V_j^{(1)}=\int_{t_n}^{t_{n+1}}V(t,x_j)\,dt\approx
\frac{\tau}{6}\left[V(t_n,x_j)+4V\left(t_n+\frac{\tau}{2},x_j\right)+V(t_{n+1},x_j)\right],\\
&A_{1,j}^{(1)}=\int_{t_n}^{t_{n+1}}A_1(t,x_j)\,dt\approx
\frac{\tau}{6}\left[A_1(t_n,x_j)+4A_1\left(t_n+\frac{\tau}{2},x_j\right)+A_1(t_{n+1},x_j)\right].
\end{align*}
\end{remark}

\subsection{Mass conservation and optimal error estimates}

Similar to the TSFP  for the Dirac equation in \cite{BCJ},
we can show that the TSFP (\ref{nleq:tsfp})  for the NLDE conserves the mass
in the discretized level with the details omitted here for brevity.

\begin{lemma}
\label{nlMass_STA_TSFP}
The  TSFP (\ref{nleq:tsfp}) conserves
the mass in the discretized level, i.e.
\be
\|\Phi^n\|_{l^2}^2:=h\sum_{j=0}^{M-1}|\Phi_j^{n}|^2\equiv  h\sum_{j=0}^{M-1}|\Phi_j^{0}|^2
=\|\Phi^0\|_{l^2}^2=h\sum_{j=0}^{M-1}|\Phi_0(x_j)|^2,\qquad n\ge0.
\ee
\end{lemma}

From Lemma \ref{nlMass_STA_TSFP}, we conclude that the TSFP (\ref{nleq:tsfp}) is unconditionally stable.
In addition, following the error estimate of the TSFP method for the
nonlinear Schr\"{o}dinger equation (NLSE) via the formal Lie calculus introduced in \cite{Lubich,BC1},
it is easy to show the following error estimate of the TSFP for the NLDE (see its proof in Appendix C).
For the simplicity of notations, we shall only consider  the NLDE with nonlinearity $\bF(\Phi)$
given in \eqref{nolcc} with $\lambda_1=0$ and time independent potential, i.e.
$V(t,x)=V(x)$ and $A_1(t,x)=A_1(x)$. In such case, the TSFP \eqref{nleq:tsfp} is the numerical method under consideration.

We make the following assumptions
on the time-independent electromagnetic potentials
 \begin{equation*}
(E)\hskip4cm \|V(x)\|_{W_p^{m_0,\infty}}+
\|A_1(x)\|_{W_p^{m_0,\infty}}\lesssim 1,\qquad m_0\ge4,\hskip6cm
\end{equation*}
and the exact solution $\Phi:=\Phi(t,x)$ of the NLDE (\ref{NLD1d})
\begin{align*}
(F)\quad \|\Phi\|_{ L^\infty([0,T]; (H_p^{m_0})^2)}\lesssim 1,\qquad
\|\partial_t\Phi\|_{L^{\infty}([0, T];(H_p^{m_0-1})^2)}\lesssim \frac{1}{\varepsilon^2},\quad
\|\partial_{tt}\Phi\|_{L^{\infty}([0,T];(H_p^{m_0-2})^2)}\lesssim\frac{1}{\varepsilon^4}.\hskip4cm
\end{align*}
Under the assumption (F), using \eqref{obser12} and \eqref{cons}, we immediately find that
the density behaves better than the wave function as
\begin{align}\label{eq:dens}
\|\partial_t\rho(t,x)\|_{L^{\infty}([0, T];H_p^{m_0-1})}\lesssim \frac{1}{\varepsilon},\quad
\|\partial_{tt}\rho(t,x)\|_{L^{\infty}([0,T];H_p^{m_0-2})}\lesssim\frac{1}{\varepsilon^2}.
\end{align}

The error estimates can be established as follows (see their proofs in Appendix C).
\begin{theorem}
\label{thm:tsfp1}
Let $\Phi^n$ be the approximation obtained from the TSFP (\ref{nleq:tsfp}) and $\lambda_1=0$ in \eqref{nolcc} with
time-independent potentials $V(x)$ and $A_1(x)$.
Assume $0<\tau\lesssim \eps^2$, under the
assumptions $(E)$ and $(F)$, there exists $h_0>0$ and $\tau_0>0$ sufficiently small
and independent of $\varepsilon$ such that, for any $0<\varepsilon\leq1$,
when $0<h\leq h_0$ and $0<\tau\leq\tau_0$, we have the error estimates
\be\label{eq:errortsfp1}
\|\Phi(t_n,x)-(I_M\Phi^n)(x)\|_{H^s}\lesssim\frac{\tau^2}{\varepsilon^4}+h^{m_0-s},\quad s=0,1;
\quad \|\Phi^n\|_{l^{\infty}}\leq 1+M_0,\quad 0\leq n\leq\frac{T}{\tau}.
\ee
\end{theorem}
The convergence result \eqref{eq:errortsfp1} can be refined if the time step is
chosen such that $\tau=2\pi\eps^2/N$ with positive integer $N$. More
precisely, we have the following improved error bound.
\begin{theorem}
\label{thm:tsfp}
Let $\Phi^n$ be the approximation obtained from the TSFP (\ref{nleq:tsfp}) and $\lambda_1=0$ in \eqref{nolcc} with
time-independent potentials $V(x)$ and $A_1(x)$.
Assume $\tau=\frac{2\pi\eps^2}{N}$ with positive integer $N$, under the
assumptions $(E)$ and $(F)$, there exist $h_0>0$ and $\frac{1}{N_0}>0$ sufficiently small
and independent of $\varepsilon$ such that, for any $0<\varepsilon\leq1$,
when $0<h\leq h_0$ and $N\ge N_0$, we have the error estimates
\be\label{eq:errortsfp}
\|\Phi(t_n,x)-(I_M\Phi^n)(x)\|_{H^s}\lesssim\frac{\tau^2}{\varepsilon^2}+h^{m_0-s}+N^{-m^*},\quad s=0,1;
\quad \|\Phi^n\|_{l^{\infty}}\leq 1+M_0,\quad 0\leq n\leq\frac{T}{\tau},
\ee
where $m^*$ is an arbitrary positive integer.
\end{theorem}

\begin{remark}
We remark that the $W_p^{m_0,\infty}$
assumption in (E) is necessary for the exact solution $\Phi(t,x)$ belonging to $(H_p^{m_0})^2$.
 In practice, as long as the solution of the NLDE is well localized such that
the periodic truncation of the potential term $(V(t,x)I_2-A_1(t,x)\sigma_1)\Phi(t,x)$ does not introduce significant aliasing error,
we still have the error estimates in the above theorem.
\end{remark}
\begin{remark}
For the NLDE \eqref{NLD1d} with general nonlinearity $\bF(\Phi)$ \eqref{nolcc} with $\lambda_1\neq0$, the proof is similar and we omit it here
for simplicity. In addition, the error estimates hold true in higher dimensions ($d=2,3$), if we choose $s=0,1,2$.
\end{remark}
From Theorem \ref{thm:tsfp}, we can find the temporal/spatial resolution capacity of the TSFP method for the NLDE
in the nonrelativistic limit regime, which is: $h=O(1)$ and $\tau=O(\vep^2)$. In fact,
for a given accuracy bound $\delta_0>0$, the $\vep$-scalability of the TSFP is:
\begin{equation}
\tau=O\left(\vep^2 \sqrt{\delta_0}\right)=O(\vep^2), \qquad  h=O\left(\delta_0^{1/m_0}\right)=O\left(1\right),
\qquad 0<\vep\ll1.
\end{equation}

Similar to the Appendix D in \cite{BCJ} for the Dirac equation,
it is straightforward to generalize the TSFP to the NLDE (\ref{SDEdd}) in 2D and
(\ref{SDEd}) in 1D, 2D and 3D and the details are omitted here for brevity.

\section{Numerical comparisons} \label{sec5}

In this section, we compare the accuracy of
different numerical methods including the CNFD, EWI-FP and TSFP methods
for solving the NLDE (\ref{SDEdd}) in terms of the mesh size $h$, time step $\tau$
and the parameter $0<\varepsilon\le 1$.
We will pay particular attention
to the $\vep$-scalability of different methods
in the nonrelativistic limit regime, i.e. $0<\varepsilon\ll 1$.

To test the accuracy, we take $d=1$ and choose the electromagnetic potentials in the NLDE (\ref{SDEdd}) as
\begin{eqnarray}
A_1(t,x) = \frac{(x+1)^2}{1+x^2},\qquad V(t,x) = \frac{1-x}{1+x^2}, \qquad x\in{\mathbb R}, \quad t\ge0,
\end{eqnarray}
and the initial value as
\begin{equation}
\phi_1(0,x) = e^{-x^2/2}, \quad\phi_2(0,x) = e^{-(x-1)^2/2}, \qquad x\in{\mathbb  R}.
\end{equation}
The problem is solved numerically on an interval $\Omega=(-16, 16)$, i.e. $a=-16$ and $b=16$,
with periodic boundary conditions on $\partial\Omega$.
The `exact' solution  $\Phi(t,x)=(\phi_1(t, x),\phi_2(t, x))^T$
is obtained numerically by using the TSFP method with a very fine mesh size and a small time step,
e.g. $h_e = 1/16$ and $\tau_e = 10^{-7}$ for comparing with the numerical
solutions obtained by EWI-FP and TSFP, and respectively $h_e = 1/4096$
for comparing with the numerical solutions obtained by the CNFD method.
Denote $\Phi^n_{h,\tau}$ as the numerical solution
obtained by a numerical method with  mesh size $h$ and time step $\tau$.
In order to quantify the convergence, we introduce
\be
e_{h,\tau}(t_n)=\|\Phi_{h,\tau}^n-\Phi(t_n,\cdot)\|_{l^2}=\sqrt{h\sum_{j=0}^{M-1}|\Phi^n_j-\Phi(t_n,x_j)|^2}.
\ee

\begin{table}[t!]
\def\temptablewidth{1\textwidth}
\vspace{-12pt}
\caption{Spatial error analysis of
the CNFD method for the NLDE (\ref{SDEdd}).}
{\rule{\temptablewidth}{1pt}}
\begin{tabular*}{\temptablewidth}{@{\extracolsep{\fill}}cccccc}
$e_{h,\tau_e}(t=2)$  & $h_0=1/8$   & $h_0/2$   &$h_0/2^2$ &  $h_0/2^3$ & $h_0/2^4$ \\
\hline
$\varepsilon_0=1$ &  8.15E-2  &  2.02E-2  & 5.00E-3  &  1.25E-3 & 3.12E-4 \\
order & - & 2.01 & 2.01 & 2.00 & 2.00\\ \hline
$\varepsilon_0/2$ &  9.29E-2  &  2.30E-2  &  5.73E-2  &  1.43E-3 & 3.58E-4 \\
order & - & 2.01 & 2.01 & 2.00 & 2.00\\ \hline
$\varepsilon_0/2^2$ &  9.91E-2 &  2.46E-2  &  6.12E-3  &  1.53E-3 &3.82E-4 \\
order & - & 2.01 & 2.01 & 2.00 & 2.00\\ \hline
$\varepsilon_0/2^3$ &  9.89E-2  &  2.47E-2  & 6.17E-3  &  1.54E-3 & 3.85E-4 \\
order & - & 2.00 & 2.00 & 2.00 & 2.00\\ \hline
$\varepsilon_0/2^4$ &  9.87E-2  &  2.48E-2  & 6.18E-3  &  1.54E-3 & 3.83E-4 \\
order & - & 1.99 & 2.00 & 2.00 & 2.01\\
\end{tabular*}
{\rule{\temptablewidth}{1pt}}
\label{table_nonlinear_CNFD}
\end{table}

\subsection{Spatial discretization errors}

Table \ref{table_nonlinear_CNFD} lists spatial errors $e_{h,\tau_e}(t=2)$ of
the CNFD method (\ref{nlcnfd1}) for different $h$ and $\vep$ with $\tau_e=10^{-6}$
such that the temporal discretization errors are negligible, and
Table \ref{table_nonlinear_TSSP} shows similar results
for the  TSFP method (\ref{nleq:tsfp}).
The spatial discretization errors of the EWI-FP method (\ref{nlAL4})-(\ref{nlAL5})
are the same as those of the  TSFP method (\ref{nleq:tsfp}) and thus they are omitted here for
brevity.

\begin{table}[htb]
\def\temptablewidth{1\textwidth}
\vspace{-12pt}
\caption{Spatial error analysis of
the TSFP method for the NLDE (\ref{SDEdd}).}
{\rule{\temptablewidth}{1pt}}
\begin{tabular*}{\temptablewidth}{@{\extracolsep{\fill}}cccccc}
$e_{h,\tau_e}(t=2)$ & $h_0$=2   & $h_0$/2   &$h_0/2^2$ &  $h_0/2^3$  &  $h_0/2^4$ \\
\hline
   $\varepsilon_0=1$  &  1.68  &  4.92E-1  &  4.78E-2  &  1.40E-4 & 2.15E-9\\
  $\varepsilon_0/2$  &  1.48  &  3.75E-1  &  1.57E-2 &  4.24E-5 & 6.60E-10\\
  $\varepsilon_0/2^2$  &  1.21  &  2.90E-1  &  4.66E-3  &  4.91E-6 & 6.45E-10\\
  $\varepsilon_0/2^3$  &  1.37  &  2.68E-1  &  2.40E-3  &  6.00E-7 & 6.34E-10\\
  $\varepsilon_0/2^4$  &  1.41  &  2.75E-1  &  1.84E-3  &  3.06E-7 & 6.13E-10\\
  $\varepsilon_0/2^5$  &  1.45  &  2.76E-1  &  1.74E-3  &  2.37E-7 & 5.98E-10\\
   \end{tabular*}
{\rule{\temptablewidth}{1pt}}
\label{table_nonlinear_TSSP}
\end{table}

From Tables \ref{table_nonlinear_CNFD}-\ref{table_nonlinear_TSSP}, we can
draw the following conclusions on spatial discretization errors for the NLDE by using different numerical methods:

\smallskip

For any fixed $\vep=\vep_0>0$, the CNFD method (and FDTD methods)
is second-order accurate, and resp.,
the EWI-FP and TSFP methods are spectrally accurate (cf. each row in Tables \ref{table_nonlinear_CNFD}-\ref{table_nonlinear_TSSP}).
For $0<\vep\le 1$, the errors are independent of $\vep$  for the  EWI-FP and TSFP methods
(cf. each column in  Table \ref{table_nonlinear_TSSP}),
 and resp., are almost independent of $\vep$ for the CNFD method (cf. each column in  Table \ref{table_nonlinear_CNFD}).
In general, for any fixed $0<\vep\le 1$ and $h>0$,
the EWI-FP and TSFP methods perform much better than the CNFD method (and FDTD methods)
in spatial discretization.

Similar to the FDTD methods for the Dirac equation \cite{BCJ},
we can observe numerically the $\vep$-dependence in the spatial discretization error,
i.e. $\frac{1}{\vep}$ in front of $h^2$, which was proven in Theorems \ref{nlthm_cnfd}.
Again, the details are omitted here for brevity.

\begin{table}[t!]
\def\temptablewidth{1\textwidth}
\vspace{-12pt}
\caption{Temporal error analysis of
the CNFD method for the NLDE (\ref{SDEdd}).}
{\rule{\temptablewidth}{1pt}}
\begin{tabular*}{\temptablewidth}{@{\extracolsep{\fill}}cccccc}
$e_{h_e,\tau}(t=2)$ & $\tau_0$=0.1  & $\tau_0$/8 & $\tau_0/8^2$ & $\tau_0/8^3$ & $\tau_0/8^4$ \\
\hline
  $\varepsilon_0=1$  &  \underline{7.13E-2}  &  9.76E-4  &  1.52E-5  & 2.38E-7  & 3.65E-9 \\
  order & - & 2.01 & 2.00 & 2.00 & 2.01\\ \hline
  $\varepsilon_0/2$  &  4.58E-1  &  \underline{7.75E-3}  &  1.21E-4  &  1.89E-6  & 2.95E-8 \\
  order & - & 1.96 & 2.00 & 2.00 & 2.00\\ \hline
  $\varepsilon_0/2^2$  &  1.74  &  2.34E-1  &  \underline{3.86E-3}  &  6.01E-5 &  9.42E-7 \\
  order & - & 0.96 & 1.97 & 2.00 & 2.00\\ \hline
  $\varepsilon_0/2^3$  &  3.13  &  5.25E-1  &  2.07E-1  &  \underline{3.49E-3} & 5.46E-5 \\
  order & - & 0.86 & 0.45 & 1.96 & 2.00\\ \hline
  $\varepsilon_0/2^4$  &  2.34  &  1.84  &  8.16E-1  &  2.04E-1 &  \underline{3.42E-3} \\
  order & - & 0.16 & 0.39 & 0.67 & 1.97\\
\end{tabular*}
{\rule{\temptablewidth}{1pt}}
\label{table_nonlinear_CNFD-T}
\end{table}

\begin{table}[t!]
\def\temptablewidth{1\textwidth}
\vspace{-12pt}
\caption{Temporal error analysis of
the EWI-FP method for the NLDE (\ref{SDEdd}).}
{\rule{\temptablewidth}{1pt}}
\begin{tabular*}{\temptablewidth}{@{\extracolsep{\fill}}ccccccc}
$e_{h_e,\tau}(t=2)$ & $\tau_0$=0.1  & $\tau_0$/4 & $\tau_0/4^2$ & $\tau_0/4^3$ & $\tau_0/4^4$ & $\tau_0/4^5$\\
\hline
  $\varepsilon_0 = 1$  &  \underline{1.62E-1}  &  8.75E-3  &  5.44E-4  &  3.40E-5  &  2.12E-6 &  1.33E-7\\
  order & - & 2.11 & 2.00 & 2.00 & 2.00 & 2.00\\ \hline
  $\varepsilon_0/2    $  &  2.02  &  \underline{2.58E-2}  &  1.59E-3  &  9.94E-5  &  6.21E-6 &  3.88E-7\\
  order & - & 3.15 & 2.01 & 2.00 & 2.00 & 2.00\\ \hline
  $\varepsilon_0/2^2  $  &  2.11  &  2.11  &  \underline{1.12E-2}  &  6.94E-4  &  4.33E-5 &  2.71E-6\\
  order & - & 0.00 & 3.78 & 2.01 & 2.00 & 2.00\\ \hline
  $\varepsilon_0/2^3  $  &  2.12  &  2.12  &  1.52E-1  &  \underline{8.88E-3}  &  5.53E-4 &  3.45E-5\\
  order & - & 0.00 & 1.90 & 2.05 & 2.00 & 2.00\\ \hline
  $\varepsilon_0/2^4  $  &  2.06  &  2.06  &  2.06  &  1.40E-1  &  \underline{8.24E-3}&  5.13E-4 \\
  order & - & 0.00 & 0.00 & 1.94 & 2.04 & 2.00\\ \hline
  $\varepsilon_0/2^5  $  &  2.09  &  2.03  &  2.03  &  2.03  &  1.36E-1  &  \underline{8.01E-3} \\
  order & - & 0.02 & 0.00 & 0.00 & 1.95 & 2.04\\
\end{tabular*}
{\rule{\temptablewidth}{1pt}}
\label{table_nonlinear_EWI-T}
\end{table}

\subsection{Temporal discretization errors}

Table \ref{table_nonlinear_CNFD-T} lists
temporal errors $e_{h_e,\tau}(t=2)$ of the CNFD method (\ref{nlcnfd1})
for different $\tau$ and $\vep$ with mesh size $h_e=1/4096$ such that spatial discretization errors
are negligible.  Tables \ref{table_nonlinear_EWI-T} and \ref{table_nonlinear_TSSP-T} show similar results
for the EWI-FP method (\ref{nlAL4})-(\ref{nlAL5}) and the TSFP method (\ref{nleq:tsfp}) for
different $\tau$ and $\vep$ with $h_e=1/16$, respectively.

\begin{table}[t!]
\def\temptablewidth{1\textwidth}
\vspace{-12pt}
\caption{Temporal error analysis of
the TSFP method for the NLDE (\ref{SDEdd}).}
{\rule{\temptablewidth}{1pt}}
\begin{tabular*}{\temptablewidth}{@{\extracolsep{\fill}}ccccccc}
$e_{h_e,\tau}(t=2)$ & $\tau_0$=0.4  & $\tau_0/4$ & $\tau_0/4^2$ & $\tau_0/4^3$
& $\tau_0/4^4$ &$\tau_0/4^5$ \\
\hline
 $\varepsilon_0=1$  & \underline{1.60E-1} & 9.56E-3 & 5.95E-4 & 3.72E-5 & 2.32E-6 & 1.46E-7  \\
 order & - & 2.03 & 2.00 & 2.00 & 2.00 & 2.00\\ \hline
 $\varepsilon_0/2$  & 8.94E-1 & \underline{3.91E-2} & 2.40E-3 & 1.50E-4 & 9.36E-6 & 5.87E-7   \\
 order & - & 2.26 & 2.01 & 2.00 & 2.00 & 2.00\\ \hline
 $\varepsilon_0/2^2$  & 2.60 & 2.18E-1 & \underline{1.06E-2} & 6.56E-4 & 4.09E-5 & 2.56E-6 \\
 order & - & 1.79 & 2.18 & 2.01 & 2.00 & 2.00\\ \hline
 $\varepsilon_0/2^3$  & 2.28 & 2.33 & 4.84E-2 & \underline{2.58E-3} & 1.60E-4 & 9.98E-6  \\
 order & - & -0.02 & 2.79 & 2.11 & 2.01 & 2.00\\ \hline
 $\varepsilon_0/2^4$  & 1.46 & 1.28 & 1.30 & 1.15E-2 & \underline{6.19E-4} & 3.84E-5   \\
 order & - & 0.10 & -0.01 & 3.41 & 2.11 & 2.01\\ \hline
 $\varepsilon_0/2^5$  & 1.53 & 3.27E-1 & 4.06E-1 & 4.13E-1 & 2.83E-3 & \underline{1.53E-4}  \\
 order & - & 1.11 & -0.16 & -0.01 & 3.59 & 2.10\\
\end{tabular*}
{\rule{\temptablewidth}{1pt}}
\label{table_nonlinear_TSSP-T}
\end{table}

From Tables \ref{table_nonlinear_CNFD-T}-\ref{table_nonlinear_TSSP-T}, we can
draw the following conclusions on temporal discretization errors for the NLDE
by using different numerical methods:

\smallskip

(i) In the $O(1)$ speed-of-light regime,
i.e. $\vep=O(1)$, all the numerical methods including
CNFD, EWI-FP and TSFP are second-order accurate (cf. the first row in  Tables \ref{table_nonlinear_CNFD-T}-\ref{table_nonlinear_TSSP-T}).
In general, the TSFP method performs much better than the CNFD (and FDTD) and EWI-FP methods  in
temporal discretization for a fixed time step (cf. Table \ref{table_nonlinear_TSSP1}).
In the non-relativistic limit  regime, i.e. $0<\vep\ll1$, for the CNFD method (and FDTD methods),
the `correct' $\vep$-scalability is $\tau=O(\vep^3)$  which verifies our theoretical results
(cf. each diagonal in Table \ref{table_nonlinear_CNFD-T}); for
the EWI-FP and TSFP methods, the `correct'
$\vep$-scalability is $\tau=O(\vep^2)$ which again confirms our theoretical results
(cf. each diagonal in Tables \ref{table_nonlinear_EWI-T}\&\ref{table_nonlinear_TSSP-T}).
In fact, for $0<\vep\le1$, one can observe clearly second-order convergence in time for
the CNFD method (and FDTD methods) only when $0<\tau\lesssim \vep^3$ (cf. upper triangles in
 Table \ref{table_nonlinear_CNFD-T}), and resp., for the
EWI-FP and TSFP methods when $0<\tau\lesssim \vep^2$ (cf. upper triangles in  Tables \ref{table_nonlinear_EWI-T}\&\ref{table_nonlinear_TSSP-T}).
In general, for any fixed $0<\vep\le 1$ and $\tau>0$,
the TSFP method performs the best, and the EWI-FP method
performs much better than the CNFD method (and FDTD methods)
in temporal discretization (cf. Tables \ref{table_nonlinear_TSSP1}\&\ref{table_nonlinear_TSSP2}).

(ii). From Table \ref{table_nonlinear_TSSP-T}, our numerical results
confirm the  error bound \eqref{eq:errortsfp}
for the TSFP method,
which is much better than \eqref{eq:errortsfp1} for the TSFP method in the nonrelativistic limit regime.

\begin{table}[t!]
\def\temptablewidth{1\textwidth}
\vspace{-12pt}
\caption{Comparison of temporal errors of
different methods for the NLDE (\ref{SDEdd}) with $\vep = 1$.}
{\rule{\temptablewidth}{1pt}}
\begin{tabular*}{\temptablewidth}{@{\extracolsep{\fill}}ccccccc}
$\vep=1$ & $\tau_0$=0.1  & $\tau_0/2$ & $\tau_0/2^2$ & $\tau_0/2^3$
& $\tau_0/2^4$ &$\tau_0/2^5$ \\
\hline
 CNFD  & 7.13E-2 & 1.82E-2 & 4.55E-3 & 1.14E-3 & 2.84E-4 & 7.11E-5  \\
 order & - & 1.97 & 2.00 & 2.00 & 2.01 & 2.00\\ \hline
 EWI-FP & 1.62E-1 & 3.56E-2 & 8.75E-3 & 2.18E-3 & 5.44E-4 & 1.36E-4   \\
 order & - & 2.19 & 2.02 & 2.00 & 2.00 & 2.00\\ \hline
 TSFP  & 9.56E-3 & 2.40E-3 & 6.56E-4 & 1.60E-4 & 3.84E-5 & 9.47E-6  \\
 order & - & 1.99 & 1.87 & 2.04 & 2.06 & 2.02\\
\end{tabular*}
{\rule{\temptablewidth}{1pt}}
\label{table_nonlinear_TSSP1}
\end{table}

\begin{table}[t!]
\def\temptablewidth{1\textwidth}
\vspace{-12pt}
\caption{Comparison of temporal errors of different numerical methods
for the NLDE (\ref{SDEdd}) under proper $\varepsilon$-scalability. }
{\rule{\temptablewidth}{1pt}}
\begin{tabular*}{\temptablewidth}{@{\extracolsep{\fill}}cccccc}
\hline
$\ba{c}
\tau=O(\varepsilon^3)\\
\ea$ &$\ba{c}
\varepsilon_0=1\\
\tau_0=0.1\\
\ea$ &$\ba{c}
\varepsilon_0/2\\
\tau_0/8\\
\ea$ &$\ba{c}
\varepsilon_0/2^2\\
\tau_0/8^2\\
\ea$ &$\ba{c}
\varepsilon_0/2^3\\
\tau_0/8^3\\
\ea$  &$\ba{c}
\varepsilon_0/2^4\\
\tau_0/8^4\\
\ea$   \\
\hline
CNFD &  7.13E-2  &  7.75E-3  & 3.86E-3  &  3.49E-3 & 3.42E-3 \\
Order in time & - & 1.07 & 0.34 & 0.05 & 0.01 \\
\hline
\hline
 $\ba{c}
\tau=O(\varepsilon^2)\\
\ea$ &$\ba{c}
\varepsilon_0=1\\
\tau_0=0.1\\
\ea$ &$\ba{c}
\varepsilon_0/2\\
\tau_0/4\\
\ea$ &$\ba{c}
\varepsilon_0/2^2\\
\tau_0/4^2\\
\ea$ &$\ba{c}
\varepsilon_0/2^3\\
\tau_0/4^3\\
\ea$   &$\ba{c}
\varepsilon_0/2^4\\
\tau_0/4^4\\
\ea$     \\
\hline
EWI-FP  &  1.62E-1  &  2.58E-2  &  1.12E-2  &  8.88E-3 & 8.24E-3  \\
Order in time & - & 1.33 & 0.60 & 0.17 & 0.05 \\  \hline
TSFP  &  9.56E-3  &  2.40E-3  &  6.56E-4  &  1.60E-4 &3.84E-5 \\
Order in time & - & 1.00 & 0.94 & 1.02 & 1.03 \\
\end{tabular*}
{\rule{\temptablewidth}{1pt}}
\label{table_nonlinear_TSSP2}
\end{table}

\subsection{Comparison for $\vep=1$ and resonance regimes}

For comparison, Table  \ref{table_nonlinear_TSSP1}
depicts temporal errors
of different numerical methods when $\vep=1$ for different $\tau$,
and Table \ref{table_nonlinear_TSSP2}
shows the $\vep$-scalability
of different methods in the nonrelativistic limit regime.

\begin{table}[t!]
\def\temptablewidth{1\textwidth}
\vspace{-12pt}
\caption{Comparison of properties of different numerical methods for solving the NLDE (\ref{SDEdd}) (or (\ref{SDEd}))
with $M$ being the number of grid points in space. }
{\rule{\temptablewidth}{1pt}}
{\small\begin{tabular*}{\temptablewidth}{@{\extracolsep{\fill}}cccc}
Method & CNFD & EWI-FP & TSFP \\
\hline
  Time symmetric &  Yes  &  No & Yes \\
  Mass conservation &  Yes  &  No & Yes \\
  Energy conservation &  Yes  & No & No \\
 Dispersion relation &  No  &  No & Yes \\
   Time transverse invariant &  No  &  No & Yes \\ \hline
  Unconditionally stable &  Yes  &  No & Yes  \\
  Explicit scheme &  No  &  Yes  & Yes  \\
  Temporal accuracy & 2nd & 2nd & 2nd \\
  Spatial accuracy & 2nd & Spectral & Spectral \\
  Memory cost & $O(M)$ & $O(M)$ & $O(M)$ \\
  Computational cost & $\gg O(M)$ & $O(M\ln M)$ & $O(M\ln M)$ \\ \hline
  $\ba{c}
  \hbox{Resolution} \\
  \hbox{when}\, 0<\varepsilon\ll1\\
  \ea$  &$\ba{c}
   h=O(\sqrt{\varepsilon})\\
   \tau=O(\varepsilon^3)\\
   \ea$  &$\ba{c}
   h=O(1)\\
   \tau=O(\varepsilon^2)\\
   \ea$  &$\ba{c}
   h=O(1)\\
   \tau=O(\varepsilon^2)\\
   \ea$   \\
   \end{tabular*}}
{\rule{\temptablewidth}{1pt}}
\label{table_nonlinear_properties}
\end{table}

\bigskip

Based on the above comparisons, in view of both temporal and spatial accuracy
 and $\vep$-scalability,  we conclude that
the TSFP and EWI-FP  methods perform much better than the CNFD method (and FDTD methods)
for the discretization of the NLDE (\ref{SDEdd}) (or (\ref{SDEd})), especially in the
nonrelativistic limit regime. For the reader's convenience,
we summarize the properties of different numerical methods for the NLDE in
Table \ref{table_nonlinear_properties}.

\subsection{Temporal errors on physical observables by TSFP}

As observed in \cite{BJP1,BJP2}, the time-splitting spectral (TSSP) method
for the NLSE performs much better for the
physical observables, e.g. density and current, than for the wave function,
in the semiclassical limit regime with respect to the scaled Planck constant $0<\vep\ll1$.
In order to see whether this is still valid for the TSFP method for
the NLDE in the nonrelativistic limit regime,
let $\rho^n=|\Phi^n_{h,\tau}|^2$, ${\bf J}^n=\frac{1}{\vep}(\Phi^n_{h,\tau})^*\sigma_1\Phi^n_{h,\tau}$
with $\Phi^n_{h,\tau}$  the numerical solution
obtained by the TSFP method with  mesh size $h$ and time step $\tau$, and define the errors
\begin{equation*}
e_{\rho}^{h,\tau}(t_n) := \|\rho^n-\rho(t_n,\cdot)\|_{l^1}=h\sum_{j=0}^{N-1}|\rho_j^n-\rho(t_n,x_j)|,
\quad  e_{{\bf J}}^{h,\tau}(t_n) := \|{\bf J}^n-{\bf J}(t_n,\cdot)\|_{l^1}=
h\sum_{j=0}^{N-1}|{\bf J}_j^n-{\bf J}(t_n,x_j)|.
\end{equation*}
Table \ref{table_linear_TSSPd} lists temporal errors $e_{\rho}^{h_e,\tau}(t=2)$ and
$e_{\bf J}^{h_e,\tau}(t=2)$ of the TSFP method (\ref{nleq:tsfp}) for
 different $\tau$ and $\vep$ with $h_e=1/16$.

 \begin{table}[t!]
\def\temptablewidth{1\textwidth}
\vspace{-12pt}
\caption{Temporal errors for density and current of the TSFP for the NLDE (\ref{SDEdd}).}
{\rule{\temptablewidth}{1pt}}
\begin{tabular*}{\temptablewidth}{@{\extracolsep{\fill}}cccccccc}
$e_{\rho}^{h_e,\tau}(t=2)$  & $\tau_0$=0.4  & $\tau_0/4$ & $\tau_0/4^2$ & $\tau_0/4^3$
& $\tau_0/4^4$ &$\tau_0/4^5$ & $\tau_0/4^6$ \\
\hline
$\varepsilon_0=1$  & \underline{2.43E-1} & 1.51E-2 & 9.40E-4 & 5.88E-5 & 3.67E-6 & 2.29E-7 & 1.42E-8 \\
order & - & 2.00 & 2.00 & 2.00 & 2.00 & 2.00 & 2.01 \\ \hline
 $\varepsilon_0/2$  & 1.08 & \underline{3.32E-2} & 2.04E-3 & 1.28E-4 & 7.98E-6 & 4.98E-7 & 3.09E-8 \\
order & - & 2.51 & 2.01 & 2.00 & 2.00 & 2.00 & 2.01\\ \hline
 $\varepsilon_0/2^2$  & 1.53 & 1.67E-1 & \underline{7.68E-2} & 4.72E-4 & 2.95E-5 & 1.84E-6 & 1.16E-7  \\
order & - & 1.60 & 0.56 & 3.67 & 2.00 & 2.00 & 1.99\\ \hline
 $\varepsilon_0/2^3$  & 1.30 & 1.27 & 2.14E-2 & \underline{9.68E-4} & 5.95E-5 & 3.72E-6 & 2.32E-7 \\
order & - &  0.02 & 2.95 & 2.23 & 2.01 & 2.00 & 1.99\\ \hline
 $\varepsilon_0/2^4$  & 1.25 & 9.44E-1 & 9.40E-1 & 5.81E-3 & \underline{2.74E-4} & 1.69E-5 & 1.06E-6 \\
order & - & 0.20 & 0.00 & 3.67 & 2.20 & 2.01 & 2.00\\ \hline
 $\varepsilon_0/2^5$  & 1.13 & 3.41E-1 & 3.27E-1 & 3.27E-1 & 1.38E-3 & \underline{6.58E-5} & 4.06E-6 \\
order & - & 0.20 & 0.86 & 0.03 & 3.94 & 2.20 & 2.01\\ \hline
\end{tabular*}
{\rule{\temptablewidth}{1pt}}
\begin{tabular*}{\temptablewidth}{@{\extracolsep{\fill}}cccccccc}
$e_{\bf J}^{h_e,\tau}(t=2)$ & $\tau_0$=0.4  & $\tau_0/4$ & $\tau_0/4^2$ & $\tau_0/4^3$
& $\tau_0/4^4$ &$\tau_0/4^5$ & $\tau_0/4^6$ \\  \hline
 $\varepsilon_0=1$  & \underline{1.23E-1} & 7.20E-3 & 4.47E-4 & 2.79E-5 & 1.74E-6 & 1.09E-7 & 6.73E-9 \\
order & - & 2.05 & 2.00 & 2.00 & 2.00 & 2.00 & 2.01\\   \hline
 $\varepsilon_0/2$  & 9.64E-1 & \underline{4.38E-2} & 2.67E-3 & 1.67E-4 & 1.04E-5 & 6.50E-7 & 4.07E-8 \\
order & - & 2.23 & 2.02 & 2.00 & 2.00 & 2.00 & 2.00\\ \hline
 $\varepsilon_0/2^2$  & 1.91 & 1.81E-1 & \underline{8.19E-3} & 5.03E-4 & 3.14E-5 & 1.96E-6 & 1.23E-7  \\
order & - & 1.70 & 2.23 & 2.01 & 2.00 & 2.00 & 2.00\\ \hline
 $\varepsilon_0/2^3$  & 1.53 & 1.59 & 3.40E-2 & \underline{1.65E-3} & 1.02E-4 & 6.37E-6 & 3.98E-7 \\
order & - & -0.03 & 2.77 & 2.18 & 2.01 & 2.00 & 2.00\\ \hline
 $\varepsilon_0/2^4$  & 1.00 & 1.43 & 1.44 & 9.42E-3 & \underline{5.04E-4} & 3.13E-5 & 1.95E-6 \\
order & - & -0.26 & -0.01 & 3.63 & 2.11 & 2.00 & 2.00\\ \hline
 $\varepsilon_0/2^5$  & 1.24 & 3.10E-1 & 3.71E-1 & 3.76E-1 & 1.59E-3 & \underline{8.48E-5} & 5.26E-6 \\
order & - & 1.00 & -0.13 & -0.01 & 3.94 & 2.11 & 2.01\\
\end{tabular*}
{\rule{\temptablewidth}{1pt}}
\label{table_linear_TSSPd}
\end{table}

 From this table, we can see that
the approximations of the density and current are at the same order as for the wave function
by using the TSFP method. The reason that we can speculate is
that $\rho=O(1)$ and ${\bf J} =O(\vep^{-1})$
(see details in (\ref{obser11}) or (\ref{obser12})) in
the NLDE, while in the NLSE both density and current
are at $O(1)$, when $0<\vep\ll1$.
Furthermore, by using the results in Theorems \ref{thm:tsfp1} and \ref{thm:tsfp},
we can immediately obtain the following error bounds for the density
under the conditions in Theorem \ref{thm:tsfp1}
\be\label{notssp3}
\|\rho(t_n,\cdot)-I_M\rho^n\|_{H^s}\lesssim \frac{\tau^2}{\varepsilon^4}+h^{m_0-s},\quad s=0,1,
\quad 0\leq n\leq\frac{T}{\tau},
\ee
and respectively, under the conditions in Theorem \ref{thm:tsfp}
\be
\|\rho(t_n,\cdot)-I_M\rho^n\|_{H^s}\lesssim \frac{\tau^2}{\varepsilon^2}+h^{m_0-s}+N^{-m^*},\quad s=0,1;
\quad \|\Phi^n\|_{l^{\infty}}\leq 1+M_0,\quad 0\leq n\leq\frac{T}{\tau}.
 \ee

\section{Conclusion}\label{sec6}
Three types of numerical methods based on different space discretizations and time integrations
were analyzed rigorously and compared numerically  for solving the nonlinear Dirac equation (NLDE)
 in the nonrelativistic limit regime, i.e.  $0<\varepsilon\ll1$.
The first class is  the second order CNFD method (and FDTD including LFFD and SIFD methods).
The error estimate of the CNFD method was
rigorously analyzed, which suggests that the $\vep$-scalability of the CNFD (and FDTD)  is
$\tau=O(\varepsilon^3)$ and  $h=O(\sqrt{\varepsilon})$. The second class applies the
Fourier spectral discretization in space and  Gautschi-type integration in time,
resulting in the EWI-FP method. Rigorous error bounds for the EWI-FP method were derived,
which show that  the $\varepsilon$-scalability of the EWI-FP method is  $\tau=O(\varepsilon^2)$
 and $h=O(1)$. The last class combines the Fourier spectral discretization in space and  time-splitting technique in time, which leads to the TSFP method. Based on the rigorous error analysis, the $\varepsilon$-scalability
 of the TSFP method is  $\tau=O(\varepsilon^2)$ and  $h=O(1)$,
 which is similar to the EWI-FP method.
 From the error analysis and numerical results, the TSFP and EWI-FP  methods perform
 much better than the CNFD (and FDTD methods), especially in the nonrelativistic limit regime.
Extensive numerical results indicate that the TSFP method is superior than the EWI-FP in terms of accuracy and efficiency, and thus the TSFP method is favorable for
solving the NLDE directly, especially in the nonrelativistic limit regime.

\bigskip

\setcounter{equation}{0}  

\begin{center}
{\bf Appendix A}. Proof of Theorem \ref{nlthm_cnfd} for the CNFD method
\end{center}
\setcounter{equation}{0}
\renewcommand{\theequation}{A.\arabic{equation}}

{\it Proof.} Compared to the proof of the CNFD method for
the Dirac equation in \cite{BCJ}, the main difficulty
is to show the numerical solution $\Phi^n$ is uniformly bounded, i.e.
$\|\Phi^n\|_{l^\infty}\lesssim 1$.
In order to do so, we adapt the cut-off technique to
truncate the nonlinearity $\bF(\Phi)$ to a
global Lipschitz function with compact support \cite{BC1,BC2,BC3}.
Choose a smooth function $\alpha(\rho) (\rho \ge 0) \in C^{\infty}([0,\infty))$ defined as
\begin{equation}
\alpha(\rho)= \left\{ \begin{array}{l}
1,\qquad 0\leq\rho\leq 1, \\
 \in [0, 1],\quad 1\leq \rho\leq 2,\\
0,\qquad \rho\geq 2.
\end{array} \right.
\end{equation}
Denote $M_1=2(1+M_0)^2>0$ and define
\begin{equation}
\label{trun_nl}
\bF_{M_1}(\Phi)=\alpha\left(\frac{|\Phi|^2}{M_1}\right)\bF(\Phi), \qquad \Phi\in {\mathbb C}^2,
\end{equation}
then $\bF_{M_1}(\Phi)$ has compact support and is smooth and global Lipschitz, i.e.,
\begin{equation}
\label{GL}
\|\bF_{M_1}(\Phi_1)-\bF_{M_1}(\Phi_2)\|\leq C_{M_1}\Bigl|\, \Phi_1-\Phi_2\, \Bigr|\lesssim
\Bigl|\, \Phi_1-\Phi_2\, \Bigr|, \qquad
\Phi_1,\Phi_2\in {\mathbb C}^2,
\end{equation}
where $C_{M_1}$ is a constant independent of $\vep$, $h$ and $\tau$.
Choose $\widetilde{\Phi}^n\in X_M$ ($n\ge 0$) such that $\widetilde{\Phi}^0=\Phi^0$ and
$\widetilde{\Phi}^n$ ($n\ge1$), with $\widetilde{\Phi}^n=(\widetilde{\Phi}^n_0,\widetilde{\Phi}^n_1,\ldots,
\widetilde{\Phi}^n_M)^T$ and $\widetilde{\Phi}^n_j=(\widetilde{\phi}_{1,j}^n,\widetilde{\phi}_{2,j}^n)^T$
for $j = 0, 1, ... M$, be the numerical solution of the following finite difference equation
\begin{align}\label{NCN1}
i\delta_t^+\tilde{\Phi}_{j}^{n}=\left[-\frac{i}{\vep}\sigma_1\delta_x +\frac{1}{\vep^2}\sigma_3
+V_j^{n+1/2}I_2-A_{1,j}^{n+1/2}\sigma_1+\bF^{n+1/2}_{M_1,j}\right]\tilde{\Phi}_{j}^{n+1/2},\quad 0\le j\le M-1,\
n\ge0,
\end{align}
where $\tilde{\Phi}_{j}^{n+1/2}=\frac{1}{2}\left[\tilde{\Phi}_{j}^{n}+\tilde{\Phi}_{j}^{n+1}\right]$ and
$\bF^{n+1/2}_{M_1,j}=\frac{1}{2}\left[\bF_{M_1}(\tilde{\Phi}_{j}^{n})+\bF_{M_1}(\tilde{\Phi}_{j}^{n+1})\right]$
for $j=0,1,\ldots,M$.
In fact, we can view $\widetilde{\Phi}^n$ as
another approximation to $\Phi(t_n,x)$. Define the corresponding errors:
\begin{eqnarray*}
\widetilde{\bee}_j^n = \Phi(t_n, x_j)-\widetilde{\Phi}_j^n,\qquad j = 0, 1, ..., M,\qquad n\geq 0
\end{eqnarray*}
Then the local truncation error $\widetilde{\xi}^n\in X_M$ of the scheme (\ref{NCN1}) is defined as
\bea
\label{NTE1}
\widetilde{\xi}_j^n:=&i\delta^+_t\Phi(t_n,x_j)-
\left[-\frac{i}{\vep}\sigma_1\delta_x +\frac{1}{\vep^2}\sigma_3
+V_j^{n+1/2}I_2-A_{1,j}^{n+1/2}\sigma_1+\bW_j^n(\Phi)\right]\frac{\Phi(t_{n+1},x_j)+\Phi(t_n,x_j)}{2}, \quad
\eea
where
\be\label{gjnt67}
\bW_j^n(\Phi)=\frac{1}{2}\left[\bF_{M_1}(\Phi(t_n,x_j))+\bF_{M_1}(\Phi(t_{n+1},x_j))\right],
\qquad j=0,1,\ldots,M, \quad n\ge0.
\ee
Taking the Taylor expansion in the local truncation error \eqref{NTE1}, noticing (\ref{NLD1d})
and (\ref{trun_nl}), under the assumptions  $(A)$ and $(B)$,
with the help of triangle inequality and Cauchy-Schwartz inequality, we have
\begin{align}\label{xijn78}
|\widetilde{\xi}_{j}^n|\le
&\frac{\tau^2}{24}\|\partial_{ttt}\Phi\|_{L^{\infty}(\overline{\Omega}_T)}+
\frac{h^2}{6\varepsilon}\|\partial_{xxx}\Phi\|_{L^{\infty}(\overline{\Omega}_T)}+
\frac{\tau^2}{8\varepsilon}\|\partial_{xtt}\Phi\|_{L^{\infty}(\overline{\Omega}_T)}\nonumber\\
&+\frac{\tau^2}{8}\left(\frac{1}{\vep^2}+2+2(|\lambda_1|+|\lambda_2|)M_0^2+V_{\rm max}+A_{1,\rm max}\right)
\|\partial_{tt}\Phi\|_{L^{\infty}(\overline{\Omega}_T)}\nonumber\\
\lesssim&\frac{\tau^2}{\varepsilon^6}+\frac{h^2}{\varepsilon}+\frac{\tau^2}{\varepsilon^5}
+\frac{\tau^2}{\varepsilon^6}\lesssim\frac{\tau^2}{\varepsilon^6}+\frac{h^2}{\varepsilon},
\qquad j=0,1,\ldots, M-1, \quad n\ge0.
\end{align}
Subtracting (\ref{NTE1}) from (\ref{NCN1}), we can obtain
\begin{align}
\label{NERR1}
i\delta^+_t\widetilde{\bee}_j^n=&\left[-\frac{i}{\vep}\sigma_1\delta_x +\frac{1}{\vep^2}\sigma_3
+V_j^{n+1/2}I_2-A_{1,j}^{n+1/2}\sigma_1\right]\tilde{\bee}_{j}^{n+1/2}
+\widetilde{\xi}_j^{n}+\tilde{\eta}_{j}^{n}, \quad 0\le j\le M-1,\quad n\ge0,
\end{align}
where $\tilde{\bee}_{j}^{n+1/2}=\frac{1}{2}\left[\tilde{\bee}_{j}^{n}+\tilde{\bee}_{j}^{n+1}\right]$ and
\be \label{etatt45}
\widetilde{\eta}_{j}^n=\frac{1}{2}\bW_j^n(\Phi)\left[\Phi(t_{n+1},x_j)+\Phi(t_n,x_j)\right]-
\bF^{n+1/2}_{M_1,j}\tilde{\Phi}_{j}^{n+1/2},\quad 0\le j\le M-1,\quad n\ge0.
\ee
Combining (\ref{etatt45}), (\ref{gjnt67}) and (\ref{GL}), we get
\be
\left|\tilde{\eta}_{j}^{n}\right|\lesssim
|\widetilde{\bee}_j^{n+1}|+|\widetilde{\bee}_j^n|, \qquad 0\le j\le M-1,\quad n\ge0.
\ee
Multiplying both sides of (\ref{NERR1}) by $h({\widetilde{\bee}_j}^{n+\frac{1}{2}})^*$,
summing them up for $j=0, 1, .., M-1$, taking imaginary parts and
applying the Cauchy inequality, noticing (\ref{xijn78}), we can have
\begin{eqnarray}
\|\widetilde{\bee}^{n+1}\|^2_{l^2}-\|\widetilde{\bee}^n\|^2_{l^2}&\lesssim&\tau
\left(\|\widetilde{\xi}^n\|_{l^2}^2+\|\widetilde{\xi}^n\|_{l^2}^2+
\|\widetilde{\bee}^{n+1}\|_{l^2}^2+\|\widetilde{\bee}^n\|_{l^2}^2\right)\nonumber\\
&\lesssim&\tau\left[\left(\frac{h^2}{\varepsilon}+\frac{\tau^2}{\varepsilon^6}\right)^2
+\|\widetilde{\bee}^{n+1}\|_{l^2}^2+\|\widetilde{\bee}^n\|_{l^2}^2\right],\qquad n\ge0.
\end{eqnarray}
Summing the above inequality, we obtain
\be
\|\widetilde{\bee}^{n}\|^2_{l^2}-\|\widetilde{\bee}^0\|^2_{l^2}
\lesssim \tau \sum_{l=0}^n\|\widetilde{\bee}^{l}\|^2_{l^2}+\left(\frac{h^2}{\varepsilon}
+\frac{\tau^2}{\varepsilon^6}\right)^2,\qquad 0\le n \le \frac{T}{\tau}.
\ee
Using the discrete Gronwall's inequality and noting $\widetilde{\bee}^0={\bf 0}$,
there exist $0<\tau_1\le \frac{1}{2}$ and $h_1>0$ sufficiently small and independent of $\vep$,
when $0<\tau\leq\tau_1$ and $0<h\le h_1$, we get
\begin{equation}
\|\widetilde{\bee}^n\|_{l^2}\lesssim \frac{h^2}{\varepsilon}+\frac{\tau^2}{\varepsilon^6},
\qquad 0\le n \le \frac{T}{\tau}.
\end{equation}
Applying the inverse inequality in 1D, we have
\begin{equation}
\|\widetilde{\bee}^n\|_{l^{\infty}}\lesssim\frac{1}{\sqrt{h}}\|\widetilde{\bee}^n\|_{l^2}\lesssim \frac{h^{\frac{3}{2}}}{\varepsilon}+\frac{\tau^2}{\varepsilon^6\sqrt{h}},\qquad 0\le n \le \frac{T}{\tau}.
\end{equation}
Under the conditions $0<\tau\lesssim \vep^3h^{\frac{1}{4}}$
and $0<h\lesssim \vep^{2/3}$, there exist $h_2>0$ and $\tau_2>0$ sufficiently small and independent of $\vep$,
when $0<h\leq h_2$ and $0<\tau\le \tau_2$, we get
\begin{eqnarray}
\|\widetilde{\Phi}^n\|_{l^{\infty}}\leq \|\Phi\|_{L^{\infty}(\Omega_T)}+\|\widetilde{\bee}^n\|_{l^{\infty}}\leq
1+M_0, \qquad 0\le n \le \frac{T}{\tau}.
\end{eqnarray}
Therefore, under the conditions in Theorem \ref{nlthm_cnfd},
the discretization (\ref{NCN1}) collapses exactly to the CNFD discretization
(\ref{nlcnfd1}) for the NLDE if we take $\tau_0=\min\{1/2,\tau_1,\,\tau_2\}$ and
$h_0=\min\{h_1,\, h_2\}$, i.e.
\begin{equation}
\widetilde{\Phi}^n = \Phi^n, \qquad 0\leq n \leq\frac{T}{\tau}.
\end{equation}
Thus the proof is completed. \hfill $ \square$

\bigskip


\bigskip

\begin{center}
{\bf Appendix B}. Proof of Theorem \ref{nlthm_EWI} for the EWI-FP method
\end{center}
\setcounter{equation}{0}
\renewcommand{\theequation}{B.\arabic{equation}}
{\it Proof.} Here the main difficulty
is to show that the numerical solution $\Phi^n_M(x)$ is uniformly bounded, i.e.
$\|\Phi^n_M(x)\|_{L^\infty}\lesssim 1$, which will be established by  the method
of mathematical induction \cite{BC1,BC2,BC3}.
Define the error function $\bee^n(x)\in Y_M$ for $n\ge0$ as
\begin{equation}
\label{nlpf_EWI_1}
\bee^n(x)=P_M\Phi(t_n,x)-\Phi_M^n(x)=
\sum_{l=-M/2}^{M/2-1}\widehat{\bee}_l^n e^{i\mu_l(x-a)},\qquad a\leq x\leq b,\quad n\ge0.
\end{equation}
Using the triangular inequality and standard interpolation result, we get
\bea
\label{nlError_L2}
\|\Phi( t_n,x)-\Phi_M^n(x)\|_{L^2}
\leq \|\Phi( t_n,x)-P_M\Phi(t_n,x)\|_{L^2}+\|\bee^n(x)\|_{L^2}
\leq h^{m_0}+\|\bee^n(x)\|_{L^2},
\ \ 0\leq n\leq\frac{T}{\tau}.\
\eea
Thus we only need estimate $\|\bee^n(x)\|_{L^2}$.
It is easy to see that (\ref{nlthm_eq_EWI}) is valid when $n=0$.

Define the local truncation error $\xi^n(x)=\sum\limits_{l=-M/2}^{M/2-1}
\widehat{\xi}_l^ne^{i\mu_l(x-a)}\in Y_M$ of the
EWI-FP (\ref{nlCoe2}) for $n\ge0$ as
\be\label{nleq:localerr}
\widehat{\xi}_l^n=\begin{cases}\widehat{(\Phi(\tau))}_l
-e^{-i\tau\Gamma_l/\vep^2}\widehat{(\Phi(0))}_l
+i\vep^2\Gamma_l^{-1}\left[I_2-e^{-\frac{i\tau}{\vep^2}\Gamma_l}\right]
\widehat{\bG(\Phi)}_l(0), &n=0,\\
\\
\widehat{(\Phi(t_{n+1}))}_l-
e^{-i\tau\Gamma_l/\vep^2}\widehat{(\Phi(t_n))}_l+iQ_l^{(1)}(\tau)
\widehat{\bG(\Phi)}_l(t_n)+iQ_l^{(2)}(\tau)\delta_t^-\widehat{\bG(\Phi)}_l(t_n), &n\ge1,
\end{cases}
\ee
where we denote $\Phi(t)$ and $\bG(\Phi)$ in short for $\Phi(t,x)$ and
$\bG(\Phi(t,x))$ in (\ref{bGphi789}), respectively,
for the simplicity of notations.
In order to estimate the local truncation error $\xi^n(x)$,
multiplying both sides of the NLDE (\ref{NLD1d})
by $e^{i\mu_l(x-a)}$ and integrating over the interval $(a,b)$,
we easily recover the equations for $\widehat{\Phi(t)}_l$,
which are exactly the same as (\ref{nlODE765}) with $\Phi_M$ being
replaced by $\Phi(t,x)$.  Replacing $\Phi_M$ with $\Phi(t,x)$, we
use the same notations $\widehat{\bG(\Phi)_l^n}(s)$  as in (\ref{nleq:fdef}) and
 the time derivatives of  $\widehat{\bG(\Phi)_l^n}(s)$
enjoy the same properties  of time derivatives of $\Phi(t,x)$.
Thus, the same representation (\ref{nlEWIn})
holds for $\widehat{\Phi(t_n)}_l$ for $n\ge1$.
From the derivation  of the EWI-FS method,
it is clear that the error $\xi^n(x)$ comes from the approximations
for the integrals in (\ref{nleq:inteapp0}) and (\ref{nleq:inteappn}).
Thus we have
\begin{align}
\label{EWI_Tr_1:1}
\widehat{\xi}^0_l=&-i\int_0^{\tau}e^{\frac{i(s-\tau)}
{\varepsilon^2}\Gamma_l}\left[\widehat{\bG(\Phi)^0_l}(s)
-\widehat{\bG(\Phi)^0_l}(0)\right]ds=-i\int_0^{\tau}\int_0^se^{\frac{i(s-\tau)}
{\varepsilon^2}\Gamma_l}\partial_{s_1}\widehat{\bG(\Phi)^0_l}(s_1)\,ds_1ds,
\end{align}
and for $n\ge1$
\begin{align}
\widehat{\xi}^n_l=&-i\int_0^{\tau}e^{\frac{i(s-\tau)}{\varepsilon^2}\Gamma_l}
\left(\int_0^s\int_0^{s_1}\partial_{s_2s_2}\widehat{\bG(\Phi)^{n}_l}(s_2)\,ds_2ds_1
+s\int_0^1\int_{\theta\tau}^\tau\partial_{\theta_1\theta_1}\widehat{\bG(\Phi)_l^{n-1}}
(\theta_1)\,d\theta_1d\theta\right)ds.\label{EWI_Tr_1:2}
\end{align}
Subtracting (\ref{nlCoe2}) from (\ref{nleq:localerr}),
we obtain
\bea
\label{nlError_fun_EWI}
&&\widehat{\bee}^{n+1}_l=
e^{-i\tau\Gamma_l/\vep^2}\widehat{\bee}^n_l+\widehat{R}^n_l
+\widehat{\xi}^n_l,\qquad 1\leq n\leq\frac{T}{\tau}-1,\\
\label{nlError_fun_EWI7}
&&\widehat{\bee}^0_l={\bf 0},\qquad\widehat{\bee}^1_l=
\widehat{\xi}^0_l,\qquad l=-\frac{M}{2},...,\frac{M}{2}-1.
\eea
where $R^n(x)=\sum\limits_{l=-M/2}^{M/2-1}\widehat{R}_l^ne^{i\mu_l(x-a)}\in Y_M$ for $n\ge1$ is given by
\be\label{nleq:R1}
\widehat{R}^n_l=-iQ_l^{(1)}(\tau)\left[\widehat{\bG(\Phi(t_n))}_l-\widehat{\bG(\Phi^n_M)}_l\right]
-iQ_l^{(2)}(\tau)\left[\delta_t^-\widehat{\bG(\Phi(t_n))}_l-\delta_t^-\widehat{\bG(\Phi^n_M)}_l\right].
\ee
From (\ref{EWI_Tr_1:1}) and (\ref{nlError_fun_EWI7}), we have
\be\label{ert987}
|\widehat{\xi}^0_l|\lesssim\int_0^{\tau}\int_0^s\left|\partial_{s_1}\widehat{\bG(\Phi)^0_l}(s_1)\right|ds_1ds.
\ee
By the Parseval equality and assumptions (C) and (D), we get
\begin{align} \label{ert988}
\|\bee^1(x)\|_{L^2}^2=&\|\xi^0(x)\|_{L^2}^2=(b-a)\sum\limits_{l=-M/2}^{M/2-1}\left|\widehat{\xi}_l^0\right|^2
\lesssim (b-a)\tau^2\int_0^\tau\int_0^s\sum\limits_{l=-M/2}^{M/2-1}\left|\partial_{s_1}
\widehat{\bG(\Phi)^0_l}(s_1)\right|^2\,ds_1ds\nn \\
\lesssim& \tau^2\int_0^\tau\int_0^s\|\partial_{s_1}(\bG(\Phi(s_1))\|_{L^2}^2\,ds_1ds
\lesssim \frac{\tau^4}{\varepsilon^4}\lesssim \frac{\tau^4}{\varepsilon^8}.
\end{align}
Thus we have
\be
\|\Phi(t_1,x)-\Phi_M^1(x)\|_{L^2}\lesssim h^{m_0}+\|\bee^1(x)\|_{L^2}\lesssim
h^{m_0}+\frac{\tau^2}{\vep^4}.
\ee
By using the inverse inequality, we get
\be
\|\bee^1(x)\|_{L^\infty}\le \frac{1}{h^{1/2}}\|\bee^1(x)\|_{L^2}\lesssim \frac{\tau^2}{\vep^4h^{1/2}},
\ee
which immediately implies
\bea
\|\Phi_M^1(x)\|_{L^\infty}&\le&\|\Phi(t_1,x)\|_{L^\infty}+\|\Phi(t_1,x)-P_M\Phi(t_1,x)\|_{L^\infty}
+\|\bee^1(x)\|_{L^\infty}\nn\\
&\le&M_0+h^{m_0-1}+\frac{\tau^2}{\vep^4h^{1/2}}.
\eea
Under the conditions in Theorem \ref{nlthm_EWI}, there exist $h_1>0$ and $\tau_1>0$ sufficiently small
and independent of $\vep$, for
$0<\vep\le 1$, when $0<h\le h_1$ and $0<\tau\le \tau_1$,
we have
\be
\|\Phi_M^1(x)\|_{L^\infty}\le 1+M_0,
\ee
thus (\ref{nlthm_eq_EWI}) is valid when $n=1$.

Now we assume that (\ref{nlthm_eq_EWI}) is valid for all
$0\leq n\leq m\le \frac{T}{\tau}-1$, then we need to show that it is still valid when $n=m+1$.
Similar to (\ref{ert987}) and (\ref{ert988}), under the assumptions (C) and (D), we obtain
\be
|\widehat{\xi}^n_l|\le\int_0^{\tau}
\left(\int_0^s\int_0^{s_1}\left|\partial_{s_2s_2}\widehat{\bG(\Phi)^{n}_l}(s_2)\right|\,ds_2ds_1
+s\int_0^1\int_{\theta\tau}^\tau\left|\partial_{\theta_1\theta_1}\widehat{\bG(\Phi)_l^{n-1}}
(\theta_1)\right|\,d\theta_1d\theta\right)ds,
\ee
\bea
\|\xi^n(x)\|_{L^2}^2&=&(b-a)\sum\limits_{l=-M/2}^{M/2-1}\left|\widehat{\xi}_l^n\right|^2
\lesssim \tau^3\int_0^{\tau}\int_0^s\int_{0}^{s_1}
\sum\limits_{l=-\frac{M}{2}}^{\frac{M}{2}-1}\left|\partial_{s_2s_2}
\widehat{\bG(\Phi)_l^n}(s_2)\right|^2\,ds_2ds_1ds\nn\\
&&+\tau^3\int_0^\tau\int_{0}^1\int_{\theta\tau}^\tau s\sum\limits_{l=-\frac{M}{2}}^{\frac{M}{2}-1}
\left|\partial_{\theta_1\theta_1}\widehat{\bG(\Phi)_l^{n-1}}(\theta_1)\right|^2\,d\theta_1\,d\theta\,ds\nn\\
&\lesssim&\tau^6\|\partial_{tt}(W(\Phi(t))\|_{L^\infty([0,T]; (L^2)^2)}^2
\lesssim\frac{\tau^6}{\varepsilon^8},\qquad n=0,1,\ldots,m.
\eea
Using the properties of the matrices $Q_l^{(1)}(\tau)$ and $Q_l^{(2)}(\tau)$, it is easy to verify that
\be\label{nleq:Qbd}
\|Q_l^{(1)}(\tau)\|_2\leq \tau,\quad \|Q_l^{(2)}(\tau)\|_2\leq \frac{\tau^2}{2},\quad l=-\frac{M}{2},\ldots,\frac{M}{2}-1.
\ee
Combining (\ref{nleq:R1}) and (\ref{nleq:Qbd}), we get
\begin{align}
&\frac{1}{b-a}\|R^n(x)\|_{L^2}^2=\sum\limits_{l=-M/2}^{M/2-1}|\widehat{R}_l^n|^2\nonumber\\
&\lesssim \tau^2\sum\limits_{l=-M/2}^{M/2-1}\Biggl[\left|\widehat{(\Phi(t_n))}_l-(\widehat{\Phi_M^n})_l\right|^2
+\left|\widehat{(\Phi(t_{n-1}))}_l-(\widehat{\Phi_M^{n-1}})_l\right|^2
+\left|\widehat{\bG(\Phi)}_l(t_n)-\widehat{\bG(\Phi_M^n)}_l\right|^2\nn\\
&\ \ +\left|\widehat{\bG(\Phi)}_l(t_{n-1})-\widehat{\bG(\Phi_M^{n-1})}_l\right|^2\Biggr]\nn\\
&\lesssim\tau^2\left[\|\Phi( t_n,x)-\Phi_M^n(x)\|_{L^2}^2+\|\Phi( t_{n-1},x)-\Phi_M^{n-1}(x)\|_{L^2}^2\right]\nn\\
&\lesssim\tau^2 h^{2m_0}
+\tau^2\|\bee^{n}(x)\|_{L^2}^2+\tau^2\|\bee^{n-1}(x)\|_{L^2}^2,\qquad n=0,1,\ldots,m.\label{nleq:rnbd}
\end{align}
Multiplying both sides of (\ref{nlError_fun_EWI}) from left by $\left(\widehat{\bee}^{n+1}_l+e^{-i\tau\Gamma_l/\vep^2}\widehat{\bee}^n_l\right)^*$, taking
the real parts and using the Cauchy inequality, we obtain
\be
\left|\widehat{\bee}^{n+1}_l\right|^2-\left|\widehat{\bee}^{n}_l\right|^2
\leq
\tau\left(\left|\widehat{\bee}^{n+1}_l\right|^2+\left|\widehat{\bee}^{n}_l\right|^2
\right)+\frac{|\widehat{R}_l^n|^2}{\tau}+\frac{|\widehat{\xi}_l^n|^2}{\tau}.
\ee
Summing the above for $l=-M/2, \ldots , M/2-1$ and then multiplying it by $(b-a)$,
using the Parseval equality,  we obtain for $n\ge1$
\begin{align}\label{nleq:rec1}
\left\|\bee^{n+1}(x)\right\|_{L^2}^2-\left\|\bee^{n}(x)\right\|_{L^2}^2\lesssim  \tau\left(\left\|\bee^{n+1}(x)\right\|_{L^2}^2+\left\|\bee^{n}(x)
\right\|_{L^2}^2\right)+\frac{1}{\tau}\left(\|R^n(x)\|^2_{L^2}
+\|\xi^n(x)\|^2_{L^2}\right).
\end{align}
Summing (\ref{nleq:rec1}) for $n=1,\ldots,m$, using (\ref{nleq:rnbd}), we derive
\be
\left\|\bee^{m+1}(x)\right\|_{L^2}^2-\left\|\bee^{1}(x)\right\|_{L^2}^2\lesssim \tau\sum\limits_{k=1}^{m+1}\left\|\bee^{k}(x)\right\|_{L^2}^2+\frac{m\tau^5}{\varepsilon^8}+
m\tau h^{2m_0},\quad 1\le m\leq\frac{T}{\tau}-1.
\ee
Noticing
$\|\bee^1(x)\|_{L^2}\lesssim\frac{\tau^2}{\varepsilon^2}\lesssim\frac{\tau^2}{\varepsilon^4}$ and
using the discrete Gronwall's inequality,
there exist $0<\tau_2\le \frac{1}{2}$ and $h_2>0$ sufficiently small and independent of $\vep$ such that,
for $0<\vep\le 1$, when $0<\tau\le \tau_2$ and $0<h\le h_2$, we get
\be\label{nleq:eq:l2e}
\left\|\bee^{m+1}(x)\right\|_{L^2}^2\lesssim  h^{2m_0}+\frac{\tau^4}
{\varepsilon^8},\qquad 1\le m\leq \frac{T}{\tau}-1.
\ee
Thus we have
\be
\|\Phi(t_{m+1},x)-\Phi_M^{m+1}(x)\|_{L^2}\lesssim h^{m_0}+\|\bee^{m+1}(x)\|_{L^2}\lesssim
h^{m_0}+\frac{\tau^2}{\vep^4}.
\ee
By using the inverse inequality, we get
\be
\|\bee^{m+1}(x)\|_{L^\infty}\le \frac{1}{h^{1/2}}\|\bee^{m+1}(x)\|_{L^2}\lesssim \frac{\tau^2}{\vep^4h^{1/2}},
\ee
which immediately implies
\bea
\|\Phi_M^{m+1}(x)\|_{L^\infty}&\le&\|\Phi(t_{m+1},x)\|_{L^\infty}+\|\Phi(t_{m+1},x)-P_M\Phi(t_{m+1},x)\|_{L^\infty}
+\|\bee^{m+1}(x)\|_{L^\infty}\nn\\
&\le&M_0+h^{m_0-1}+\frac{\tau^2}{\vep^4h^{1/2}}.
\eea
Under the conditions in Theorem \ref{nlthm_EWI}, there exist $h_3>0$ and $\tau_3>0$ sufficiently small
and independent of $\vep$, for
$0<\vep\le 1$, when $0<h\le h_3$ and $0<\tau\le \tau_3$,
we have
\be
\|\Phi_M^{m+1}(x)\|_{L^\infty}\le 1+M_0,
\ee
thus (\ref{nlthm_eq_EWI}) is valid when $n=m+1$.
Then the proof of (\ref{nlthm_eq_EWI}) is completed by the method of mathematical induction under the choice of $h_0=\min\{h_1, h_2, h_3\}$ and $\tau_0=\min\{1/2, \tau_1,\tau_2, \tau_3\}$.
\hfill $\Box$

\bigskip

\begin{center}
{\bf Appendix C}. Proof of Theorems \ref{thm:tsfp} and \ref{thm:tsfp1} for the TSFP method
\end{center}
\setcounter{equation}{0}
\renewcommand{\theequation}{C.\arabic{equation}}
{\it Proof.}
As the proof of Theorem \ref{thm:tsfp} implies the conclusion of Theorem \ref{thm:tsfp1},
 we only present here the proof of Theorem \ref{thm:tsfp}  and
omit the arguments for Theorem \ref{thm:tsfp1} for brevity.

Denote ${\mathbf{T}}=(-\varepsilon\sigma_1i\partial_x+\sigma_3)/\eps^2$, where $i{\mathbf{T}}$ generates a unitary group in $(H^k_p(\Omega))^2$ ($k\ge0$). Let $\Phi^{[n]}:=\Phi^{[n]}(x)$ be the numerical approximation of $\Phi(t_n,x)$ with  $\Phi^{[0]}=\Phi_0(x)$,
\be\label{eq:tssptd}
\Phi^{<1>}=e^{-i\frac{\tau{\mathbf{T}}}{2}}\Phi^{[n]},\quad \Phi^{<2>}=e^{-i\tau\left[V(x)I_2-A_1(x)\sigma_1+\lambda_2|\Phi^{<1>}(x)|^2I_2\right]}\,\Phi^{<1>},\quad \Phi^{[n+1]}=e^{-i\frac{\tau{\mathbf{T}}}{2}}\Phi^{<2>},
\ee
where we use  $\bF(\Phi)=\lambda_2|\Phi|^2I_2$ and electro-magnetic potentials are time-independent. We can view \eqref{eq:tssptd} as a semi-discretization in time for the NLDE \eqref{NLD1d}, and respectively, the  TSFP \eqref{nleq:tsfp} as a full discretization. The proof will be divided into two parts:
(I) to prove the convergence for the above semi-discretization, and (II) to
complete the error analysis by comparing the above semi-discretization \eqref{eq:tssptd}
with the TSFP \eqref{nleq:tsfp}.

{\bf Part I} (convergence of the semi-discretization).  Under the condition that $\tau=\frac{2\pi\eps^2}{N}$ with $N$ being a positive integer, we want to show that
\begin{equation}\label{eq:tssptd1}
\|\Phi^{[n]}(\cdot)-\Phi(t_n,\cdot)\|_{H^1}\lesssim \frac{\tau^2}{\varepsilon^2}+N^{-m_\ast}, \qquad 0\leq n\leq \frac{T}{\tau},
\end{equation}
where $m_\ast$ is an arbitrary positive integer.

It is easy to verify that
 $\|\Phi^{[n]}\|_{H_p^{m_0}}\leq C_T\|\Phi_0\|_{H_p^{m_0}}$ for $m_0\ge2$ under the assumption (E).
We denote the flow $\Phi^{[n]}(x)\to \Phi^{[n+1]}(x)$ as
\be
\Phi^{[n+1]}(x)={\cal S}_\tau\left(\Phi^{[n]}\right), \qquad 0\leq n\leq \frac{T}{\tau},
\ee
and the exact solution flow $\Phi(t_n)\to \Phi(t_{n+1})$ as
\be
\Phi(t_{n+1})={\cal S}_{e,\tau}\left(\Phi(t_{n})\right)\qquad 0\leq n\leq \frac{T}{\tau}.
\ee

In 1D, it is easy to establish the following stability results in view of the Sobolev inequality $\|f\|_{\infty}^2\leq\|f\|_{L^2}\|f\|_{H^1}$ and
the fact that ${\cal S}_\tau$ and ${\cal S}_{e,\tau}$ preserve the $L^2$-norm:
\be\label{eq:stab}
\|{\cal S}_\tau(\Phi_1)-{\cal S}_\tau (\Phi_2)\|_{H^1}\leq e^{\tau C_M}\|\Phi_1-\Phi_2\|_{H^1},\quad
\|{\cal S}_{e,\tau}(\Phi_1)-{\cal S}_{e,\tau} (\Phi_2)\|_{H^1}\leq e^{\tau C_M}\|\Phi_1-\Phi_2\|_{H^1},
\ee
where $ \|\Phi_1\|_{H^1},\|\Phi_2\|_{H^1}\leq M$, and $C_M$ depend on $M$, $\|V(\cdot)\|_{W^{1,\infty}}$ and $\|A_1(\cdot)\|_{W^{1,\infty}}$.

To prove (\ref{eq:tssptd1}), we adopt the approach via formal Lie calculus introduced in \cite{Lubich}
and split the proof into three steps.

{\bf Step 1} (bounds for local truncation error).
We start with the local error, i.e. to examine the error generated by one time step evolution computed via \eqref{eq:tssptd}. Denote
\be\label{eq:v0}
\bV_0=V(x)I_2-A_1(x)\sigma_1+\lambda_2|e^{-i\frac{{\mathbf{T}}\tau}{2}}\Phi_0|^2I_2, \qquad a\le x\le b,
\ee
and
\be\label{eq:v}
\bV(s)=V(x)I_2-A_1(x)\sigma_1+\lambda_2|\Phi(s,x)|^2I_2, \qquad a\le x\le b,\quad 0\le s\le \tau.
\ee
It is clear that if $m_0\ge3$,  $\|\bV_0\|_{W^{1,\infty}}+\|\bV\|_{L^\infty([0,\tau];W^{1,\infty})}\lesssim 1$ by Sobolev embedding (the norms for $\bV_0$ and $\bV(s)$ are understood for the matrix functions).
We claim that \eqref{eq:v0} is a second order  approximation of $\bV(\tau/2)$. By Duhamel's principle and Taylor expansion, it is easy to check
\begin{align}
\Phi(\tau/2)=&e^{-i\tau{\mathbf{T}}/2}\Phi_0-i\int_0^{\frac{\tau}{2}}e^{-i(\tau/2-s){\mathrm T}}\bV(s)\Phi(s)\,ds\nn\\
=&e^{-i\tau{\mathbf{T}}/2}\Phi_0-\frac{i\tau}{2}e^{-i\tau{\mathrm T}/2}\bV(0)\Phi_0+O(\partial_{s}(\bV(s)\Phi(s))\tau^2).
\end{align}
By direct computation,  the unitary group $e^{-is{\mathbf{T}}}$ preserves the orthogonality, i.e. $\partial_s\text{Re}((e^{-is{\mathbf{T}}}\Psi_1)^*(e^{-is{\mathbf{T}}})\Psi_2)=0$.
Since $\bV(0)$ is Hermitian, we know $\Phi_0$ is orthogonal to $i\bV(0)\Phi_0$ and hence $e^{-i\tau{\mathbf{T}}/2}\Phi_0$ is orthogonal to
$-\frac{i\tau}{2}e^{-i\tau{\mathrm T}/2}\bV(0)\Phi_0$, which together with assumption (F) would imply
\begin{align}
\left\||\Phi(\tau/2)|^2-|e^{-i\tau{\mathbf{T}}/2}\Phi_0|^2\right\|_{W^{1,\infty}}=&O(\partial_{s}(\bV(s)\Phi(s))\tau^2)
\leq\tau^2\|\bV(\cdot)\Phi(\cdot)\|_{W^{1,\infty}([0,\tau/2];(H^2)^2)}\lesssim\frac{\tau^2}{\eps^2},\label{eq:densityerror}
\end{align}
which gives the second order accuracy as
\be
\|\bV_0-\bV(\tau/2)\|_{W^{1,\infty}}\lesssim \frac{\tau^2}{\eps^2}.
\ee

Next, using Taylor expansion for $e^{-i\tau\bV_0}$, we have
 \be
 {\cal S}_\tau\left(\Phi_0\right)=e^{-i{\mathbf{T}}\tau}\Phi_0-i\tau e^{-i\frac{{\mathbf{T}}\tau}{2}}\bV_0 e^{-i\frac{{\mathbf{T}}\tau}{2}}\Phi_0
 -\tau^2\int_0^1(1-\theta)e^{-i\frac{{\mathbf{T}}\tau}{2}}e^{-i\theta\tau\bV_0}\bV_0^2
 e^{-i\frac{{\mathbf{T}}\tau}{2}}\Phi_0\,d\theta.
 \ee
 On the other hand, by repeatedly using Duhamel's principle (variation-of-constant formula), we write
 \begin{align*}
 \Phi(\tau)=&e^{-i\tau{\mathbf{T}}}\Phi_0-i\int_0^\tau e^{-i(\tau-s){\mathbf{T}}}\bV(s) e^{-is {\mathbf{T}}}\Phi_0\,ds
 -\int_0^\tau\int_0^s e^{-i(\tau-s){\mathbf{T}}}\bV(s) \left(e^{-i(s-w){\mathbf{T}}}\bV(w) \Phi(w)\right)\,dw\,ds.
 \end{align*}
 Denote
 \be
 G(s)=e^{-i(\tau-s){\mathbf{T}}}\bV(s) e^{-is {\mathbf{T}}}\Phi_0,\quad B(s,w)=e^{-i(\tau-s){\mathbf{T}}}\bV(s) e^{-i(s-w){\mathbf{T}}}\bV(w) e^{-iw{\mathbf{T}}}\Phi_0,
 \ee
then the local error can be written as
\begin{align*}
{\cal S}_\tau\left(\Phi_0\right)-\Phi(\tau)=&-i\tau G\left(\frac{\tau}{2}\right)+i\int_0^\tau G(s)\,ds-\frac{\tau^2}{2}B\left(\frac{\tau}{2},\frac{\tau}{2}\right)
+\int_0^\tau\int_0^s B(s,w)\,ds\,dw+r_1+r_2+r_3,
\end{align*}
with
\begin{align*}
r_1=&-\tau^2\int_0^1(1-\theta)e^{-i\frac{\mathbf{T}\tau}{2}}e^{-i\theta\tau\bV_0}
\bV_0^2e^{-i\frac{{\mathbf{T}}\tau}{2}}\Phi_0\,d\theta+\frac{\tau^2}{2}B\left(\frac{\tau}{2},\frac{\tau}{2}\right),\\
r_2=&-\int_0^\tau\int_0^s \left(e^{-i(\tau-s){\mathbf{T}}}\bV(s) \left(e^{-i(s-w){\mathbf{T}}}\bV(w)\Phi(w)\right)-B(s,w)\right)\,dwds,\\
r_3=&-i\tau e^{-i\frac{{\mathbf{T}}\tau}{2}}\left(|e^{-i\tau{\mathbf{T}}/2}\Phi_0|^2-|\Phi(\tau/2)|^2\right)e^{-i\frac{{\mathbf{T}}\tau}{2}}\Phi_0.
\end{align*}
Using \eqref{eq:densityerror}, it is easy to verify that
\be
\|r_3\|_{H^1}\lesssim\tau^3\|\bV(\cdot)\Phi(\cdot)\|_{W^{1,\infty}([0,\tau/2];(H^2)^2)}\,\|\Phi_0\|_{H^1}\lesssim\frac{\tau^3}{\eps^2}.
\ee
Similarly, we can estimate
\begin{align}
\|r_1\|_{H^1}\lesssim &\tau^2\left\|\bV(\tau/2)-\bV_0\right\|_{W^{1,\infty}}\left(\|\bV(\tau/2)\|_{W^{1,\infty}}+
\|\bV_0\|_{W^{1,\infty}}\right)\,\|\Phi_0\|_{H^1}\nn\\
&+\tau^2\max_{\theta\in(0,1)}\left\{\left\|\partial_{\theta\theta}((1-\theta)e^{-i\frac{\mathbf{T}\tau}{2}}
e^{-i\theta\tau\bV_0}\bV_0^2e^{-i\frac{{\mathbf{T}}\tau}{2}}\Phi_0)\right\|_{H^1}\right\}\nn\\
\lesssim&\tau^2\|\bV(\cdot)\Phi(\cdot)\|_{W^{1,\infty}([0,\tau/2];(H^3)^2)}+\tau^3\left(\|\bV_0^3\Phi_0\|_{H^1}+\|\bV_0^4\Phi_0\|_{H^1}\right)
\lesssim\frac{\tau^3}{\eps^2}.
\end{align}
The quadrature rule implies
\bea
\|r_2\|_{H^1}&\lesssim&\tau^3\max\limits_{s,w\in(0,\tau)}\left\{\|e^{-i(\tau-s){\mathbf{T}}}\bV(s) (e^{-i(s-w){\mathbf{T}}}\bV(w)\partial_w\Phi(w))\|_{H^1}\right\}\nn\\
&\lesssim&\tau^3\|\partial_s\Phi(\cdot)\|_{L^\infty([0,\tau];(H^1)^2)}\lesssim\frac{\tau^2}{\eps^2},
\eea
and
\be
\left\|-\frac{\tau^2}{2}B\left(\frac{\tau}{2},\frac{\tau}{2}\right)
+\int_0^\tau\int_0^s B(s,w)\,dw\,ds\right\|_{H^1}\lesssim \tau^3\max_{0\leq w\leq s\leq\tau}\left(\|\partial_{s}B\|_{H^1}+\|\partial_{w}B\|_{H^1}\right)\lesssim
\frac{\tau^3}{\varepsilon^2},
\ee
where we use the properties of ${\mathbf{T}}$ and $B(s,w)$. We check that the above errors have the desired bounds in Theorem \ref{thm:tsfp}.
 Finally, we estimate the last term as (cf. \cite{Lubich}), which contains the major part of the local error
\be
-i\tau G(\tau/2)+i\int_0^\tau G(s)\,ds=-i\tau^3\int_0^1\text{ker}(\theta)G^{\prime\prime}(\theta\tau)\,d\theta,
\ee
and $\text{ker}(\theta)$ is the Peano kernel for midpoint rule. In addition, we have
\bea
G^{\prime\prime}(s)=-e^{-i(\tau-s){\mathbf{T}}}[{\mathbf{T}},[{\mathbf{T}},\bV(s)]] e^{-is{\mathbf{T}}}\Phi_0
+2ie^{-i(\tau-s){\mathbf{T}}}[{\mathbf{T}},\bV^\prime(s)] e^{-is{\mathbf{T}}}\Phi_0
+e^{-i(\tau-s){\mathbf{T}}}\bV^{\prime\prime}(s) e^{-is{\mathbf{T}}}\Phi_0, \quad
\eea
and the commutators can be bounded as
\begin{align}
\|[{\mathbf{T}},\bV^\prime(s)]\Psi\|_{H^1}
=&\left\|\frac{1}{\eps^2}[-\varepsilon\sigma_1i\partial_x+\sigma_3,\lambda_2\partial_s|\Phi(s)|^2 I_2]\Psi\right\|_{H^1}
=\left\|\frac{1}{\eps}[-\varepsilon\sigma_1i\partial_x,\lambda_2\partial_s|\Phi(s)|^2 I_2]\Psi\right\|_{H^1}\nn\\
\lesssim&\frac{1}{\eps}\left\|\partial_s|\Phi(\cdot)|^2\right
\|_{L^\infty([0,\tau];(W^{2,\infty})^2)}\|\Psi\|_{H^1}\lesssim\frac{1}{\eps^2}\|\Psi\|_{H^1},\\
\|\bV^{\prime\prime}(s)]\Psi\|_{H^1}\lesssim&\left\|\partial_{ss}|\Phi(\cdot)|^2\right
\|_{L^\infty([0,\tau];(W^{2,\infty})^2)}\|\Psi\|_{H^1}\lesssim\frac{1}{\eps^2}\|\Psi\|_{H^1}.
\end{align}
Here we use the properties for the density $\rho=|\Phi|^2$ in \eqref{eq:dens}. The above two estimates will yield error bounds of the type $\tau^3/\eps^2$ and then we identify that the major
obstacle in obtaining error bounds like \eqref{eq:tssptd1} is from the double commutator $[{\mathbf{T}},[{\mathbf{T}},\bV(s)]]$. Noticing that
\be
[{\mathbf{T}},[{\mathbf{T}},\bV(s)]]=[{\mathbf{T}},[{\mathbf{T}},-A_1(x)\sigma_1]]
+[{\mathbf{T}},[{\mathbf{T}},V(x)I_2]]+[{\mathbf{T}},[{\mathbf{T}},\lambda_2|\Phi(s)|^2I_2]],
\ee
since $I_2$ commutes with $\sigma_3$, by direct computation, we have
\begin{align}
&\left\|[{\mathbf{T}},[{\mathbf{T}},V(x)I_2]]\Psi\right\|_{H^1}\lesssim\frac{1}{\eps^2}
\|V(\cdot)\|_{W^{3,\infty}}\|\Psi\|_{H^1}\lesssim\frac{1}{\eps^2}\|\Psi\|_{H^1},\\
&\left\|[{\mathbf{T}},[{\mathbf{T}},\lambda_2|\Phi(s)|^2I_2]]\Psi\right\|_{H^1}\lesssim \frac{1}{\eps^2}\left\||\Phi(\cdot)|^2\right\|_{L^\infty([0,\tau];(W^{3,\infty})^2)}\|\Psi\|_{H^1}\lesssim
\frac{1}{\eps^2}\|\Psi\|_{H^1}.
\end{align}
Therefore, we further reduce the major error to the commutator $[{\mathbf{T}},[{\mathbf{T}},-A_1(x)\sigma_1]]$.
On the other hand, ${\mathrm T}$ can be expanded in phase space for each Fourier mode as
\begin{equation}
\left({\mathrm T}\Psi\right)(x)=\sum\limits_{l\in{\mathbb Z}} \frac{1}{\eps^2}\Gamma_l \widehat{\Psi}_le^{i\mu_l(x-a)},\quad a<x<b,
\end{equation}
where $\Gamma_l=Q_lD_lQ_l^*$ is given in  \eqref{nleq:Gamma} and $D_l$ is given by
\begin{equation}
 D_l=\begin{pmatrix}
\delta_l &0\\
0 &-\delta_l\\
\end{pmatrix}
=\sigma_3+\eps^2\begin{pmatrix}
\frac{\mu_l^2}{\delta_l+1} &0\\
0 &-\frac{\mu_l^2}{\delta_l+1}
\end{pmatrix}.
\end{equation}
Therefore, we have the decomposition
\be\label{eq:Tdec}
{\mathrm T}={\mathrm T}_1+{\mathrm T}_2,
\ee
where
\begin{equation}
{\mathrm T}_1\Psi(x)=\sum\limits_{l\in{\mathbb Z}} \frac{1}{\eps^2}Q_l\sigma_3Q_l^*\widehat{\Psi}_le^{i\mu_l(x-a)},\quad a<x<b;\quad {\mathrm T}_2={\mathrm T}-{\mathrm T}_1;
\end{equation}
and it is then clear that
\be
\|{\mathrm T}_2\Psi\|_{H^{s-2}}\lesssim \|\Psi\|_{H^{s}},\quad \forall\, \Psi\in H^s_p(\Omega),\quad s\ge2.
\ee
Now we write the double commutator as
\be
[{\mathbf{T}},[{\mathbf{T}},-A_1(x)\sigma_1]]\Psi=[{\mathbf{T}_1},[{\mathbf{T}_1},-A_1(x)\sigma_1]]\Psi+{\text{C}_R}\Psi,
\ee
where the residual $\text{C}_R$ part and the leading part satisfy
\be\label{eq:decbd}
\|\text{C}_R\Psi\|_{H^1}\lesssim \frac{1}{\eps^2}\|\Psi\|_{H^5},\quad\left\|[{\mathbf{T}_1},
[{\mathbf{T}_1},-A_1(x)\sigma_1]]\Psi\right\|_{H^1}\lesssim\frac{1}{\eps^4}\|\Psi\|_{H^1}.
\ee
Thus, we identify the leading error term is from $[{\mathbf{T}_1},[{\mathbf{T}_1},-A_1(x)\sigma_1]]$.
Combing all the results above, we find the one step local error as
\be\label{eq:firststep1}
{\cal S}_\tau\left(\Phi_0\right)-{\cal S}_{e,\tau}\left(\Phi_0\right)=\tilde{\Lambda}_\tau\Phi_0+r_0,
\ee
where $\|r_0\|_{H^1}\lesssim \tau^3/\eps^2$ and
\be
\tilde{\Lambda}_\tau\Phi_0=-i\tau^3\int_0^1\text{ker}(\theta)e^{-i(\tau-\theta\tau){\mathbf{T}}}[{\mathbf{T}_1},[{\mathbf{T}_1},-A_1(x)\sigma_1]] e^{-i\theta\tau{\mathbf{T}}}\Phi_0\,d\theta.
\ee
Taking the decomposition \eqref{eq:Tdec} into account, we find
\be
e^{-i(\tau-\theta\tau){\mathbf{T}}}=e^{-i(\tau-\theta\tau){\mathbf{T}_1}}+O(\mu_l^2\tau),\quad  e^{-i\theta\tau{\mathbf{T}}}=e^{-i\theta\tau{\mathbf{T}_1}}+O(\mu_l^2\tau),
\ee
which can simplify the equation \eqref{eq:firststep1} in view of $\tau\lesssim\eps^2$ and the regularity of the solution,
\be\label{eq:firststep}
{\cal S}_\tau\left(\Phi_0\right)-{\cal S}_{e,\tau}\left(\Phi_0\right)=\Lambda_\tau\Phi_0+\tilde{r}_1,
\ee
with $\|\tilde{r}_1\|_{H^1}\lesssim \tau^3/\eps^2$ and
\be
\Lambda_\tau\Phi_0=-i\tau^3\int_0^1\text{ker}(\theta)e^{-i(\tau-\theta\tau){\mathrm{T_1}}}[{\mathbf{T}_1},[{\mathbf{T}_1},-A_1(x)\sigma_1]] e^{-i\theta\tau{\mathbf{T}_1}}\Phi_0\,d\theta.
\ee
Defining ${\cal F}_s(\Psi)$ for $\Psi\in H_p^{m_0}$ and $s\in\Bbb R$ as
\be\label{eq:calFdef}
{\cal F}_{s}(\Psi)=-ie^{is{\mathbf{T}_1}}(-A_1(x)\sigma_1)e^{-is{\mathbf{T}_1}}\Psi,
\ee
and it is easy to see that ${\cal F}_{s}(\Psi)$ is a $2\pi\eps^2$ periodic function. We notice that the following also holds
\be\label{eq:Fmid}
\Lambda_\tau\Phi_0=e^{-i\tau\mathbf{T}_1}\left(\tau {\cal F}_{\tau/2}(\Phi_0)-\int_0^\tau {\cal F}_s(\Phi_0)\,ds\right).
\ee
Define the local error at $t_n$   as
\be
{\cal E}_n(x)={\cal S}_\tau(\Phi(t_{n-1},x))-\Phi(t_{n},x),\qquad a\le x\le b, \quad 1\leq n\leq\frac{T}{\tau}.
\ee
Following the above computation and \eqref{eq:decbd}, it is easy to find that
\be\label{eq:localerr}
{\cal E}_n(x)=\Lambda_\tau\Phi(t_{n-1})+\tilde{r}_n, \qquad a\le x\le b, \quad 1\leq n\leq\frac{T}{\tau},
\ee
where
\be\label{eq:lambdabd}
\|\Lambda_\tau\Phi(t_{n-1})\|_{H^1}\lesssim \frac{1}{\eps^4}\|\Phi(t_{n-1})\|_{H^1},
\quad \|\tilde{r}_n\|_{H^1}\lesssim \frac{\tau^3}{\eps^2}, \qquad 1\leq n\leq\frac{T}{\tau}.
\ee

{\bf Step 2} (bounds for the global error in one period). We study the global error for $1\leq n\leq N$ with $\tau=\frac{2\pi\eps^2}{N}$. As noticed in the above local error representation
\eqref{eq:localerr}, the leading term is \eqref{eq:Fmid}, which comes
from ${\cal F}_s(\Phi_0)$ -- a $2\pi\eps^2$ periodic function.  The problem is well suited in
such period, which is similar to the NLSE case \cite{florian}.

Under the  assumptions (E) and (F), it is easy to verify that there exists a constant $M_1>0$ independent of $\eps$  such that
\be\label{eq:bd}
\left\|{\cal S_{\tau}}^{n-k}\left({\cal S}_{e,\tau}^{n-k}(\Phi_0)\right)\right\|_{H^1}\leq M_1,\quad
\left\|{\cal S}_{e,\tau}^{n-k}\left({\cal S}_{\tau}^{n-k}(\Phi_0)\right)\right\|_{H^1}\leq M_1,\quad \forall\, 0\leq k\leq n\leq\frac{T}{\tau}.
\ee
According to the above $H^1$ bounds, we denote $e^{\tau C_1}$ as the corresponding stability constant in \eqref{eq:stab}.

Now we want to use $e^{-in\tau{\mathbf{T}}}$ as a suitable approximation of the
flow ${\cal S}_\tau^n$. To this purpose, we introduce the difference between the two flows
for  $\Psi\in H_p^{k}$ ($k\in\Bbb N$) as
\be
{\cal A}_n\Psi={\cal S}_{\tau}^n\Psi-e^{-in\tau{\mathbf{T}}}\Psi,\qquad 1\le n\leq N.
\ee
As $e^{-i\tau{\mathbf{T}}}$ preserves the $H_p^k$ norm, it is easy to find the stability
\be
\|{\cal A}_1\Psi-{\cal A}_1\Phi\|_{H^1}\leq \tau e^{\tau C_1}\|\Psi-\Phi\|_{H^1},\quad \forall \|\Psi\|_{H^1}\leq M_1, \quad \|\Phi\|_{H^1}\leq M_1,
\ee
where the constants are the same as those in \eqref{eq:bd}.

For $n\ge1$,  we use the telescope identity to obtain
\begin{align*}
{\cal S}_\tau^n\Phi-e^{-in\tau{\mathbf{T}}}\Phi=\sum\limits_{k=1}^n
\left(e^{-i\tau(n-k){\mathbf{T}}}\circ({\cal S}_{\tau}-e^{-i\tau\mathbf{T}})\circ{\cal S}_\tau^{k-1}\right)\Phi,
\end{align*}
which implies
\begin{align*}
\|{\cal A}_n\Psi-{\cal A}_n\Phi\|_{H^1}\leq&\sum\limits_{k=1}^n\left\| {\cal A}_1\left({\cal S}_\tau^{k-1}\Psi\right)-
{\cal A}_1\left({\cal S}_\tau^{k-1}\Phi\right)\right\|_{H^1}
\leq\sum\limits_{k=1}^n\tau e^{\tau C_1}\|{\cal S}_\tau^{k-1}\Psi-{\cal S}_\tau^{k-1}\Phi\|_{H^1}\\
\leq&\sum\limits_{k=1}^n\tau e^{k\tau C_1}\|\Psi-\Phi\|_{H^1}
\leq\eps^2 e^{n\tau C_1}\|\Psi-\Phi\|_{H^1},\quad 1\leq n\leq N.
\end{align*}
It is convenient to use the telescope identity to obtain the  error
\begin{align}
{\cal S}_{\tau}^{n}\Phi_0-{\cal S}_{e,\tau}^{n}\Phi_0=&\sum\limits_{k=1}^n\left(
{\cal S}_{\tau}^{k-1}\circ\left({\cal S}_\tau\circ{\cal S}_{e,\tau}^{n-k}\right)\Phi_0-{\cal S}_{\tau}^{k-1}\circ\left({\cal S}_{e,\tau}^{n-k+1}\right)\Phi_0\right)\nonumber\\
=&\sum\limits_{k=1}^n\left(
e^{-i(k-1)\tau{\mathbf{T}}}\circ{\cal E}^{n-k}\right)+\tilde{R}_n,\qquad 1\leq n\leq N, \label{eq:global}
\end{align}
where
\be
\tilde{R}_n=\sum\limits_{k=1}^n\left(
{\cal A}_{k-1}\circ{\cal S}_\tau\Phi(t_{n-k})-{\cal A}_{k-1}\circ\Phi(t_{n-k+1})\right),\qquad 1\leq n\leq N.
\ee
Using the local error bounds \eqref{eq:localerr}, we have
\begin{align}\label{eq:errod3}
\|\tilde{R}_n\|_{H^1}\leq &\sum\limits_{k=1}^n\left\|
{\cal A}_{k-1}\circ{\cal S}_\tau\Phi(t_{n-k})-{\cal A}_{k-1}\circ\Phi(t_{n-k+1})\right\|_{H^1} \nonumber \\
\leq & \sum\limits_{k=1}^n\eps^2e^{(k-1)\tau C_1}\|{\cal E}^{n-k}\|_{H^1}
\lesssim n\eps^2e^{n\tau C_1} \frac{\tau^3}{\eps^4}\lesssim e^{n\tau C_1}\tau^2.
\end{align}
Using the local error representation \eqref{eq:localerr}, we write the other term in \eqref{eq:global} as
\be\label{eq:errord1}
\sum\limits_{k=1}^n\left(
e^{-i(k-1)\tau{\mathbf{T}}}\circ{\cal E}^{n-k}\right)=\sum\limits_{k=1}^n
e^{-i(k-1)\tau{\mathbf{T}}}\Lambda_\tau \Phi(t_{n-k})+\sum\limits_{k=1}^n
e^{-i(k-1)\tau{\mathbf{T}}}\tilde{r}_{n-k+1},
\ee
and
\be\label{eq:errord2}
\left\|\sum\limits_{k=1}^n
e^{-i(k-1)\tau{\mathbf{T}}}\tilde{r}_{n-k+1}\right\|_{H^1}\leq \sum\limits_{k=1}^n
\|\tilde{r}_{n-k+1}\|_{H^1}\lesssim n\frac{\tau^3}{\eps^2}\lesssim \tau^2.
\ee
Using Duhamel's formula, we have
\be
\Phi(t_{n-k})=e^{-it_{n-k}\mathbf{T}}\Phi_0(x)-i\int_0^{t_{n-k}}(V(x)I_2-A_1(x)\sigma_1+\lambda_2|\Phi(s)|^2I_2)\Phi(s)\,ds,
\ee
where we can obtain
\begin{align}
\sum\limits_{k=1}^n
e^{-i(k-1)\tau{\mathbf{T}}}\Lambda_\tau \Phi(t_{n-k})=\sum\limits_{k=1}^ne^{-i(k-1)\tau{\mathbf{T}}}\Lambda_\tau e^{-it_{n-k}{\mathbf{T}}} \Phi_0+
\sum\limits_{k=1}^n e^{-i(k-1)\tau{\mathbf{T}}}\Lambda_\tau\hat{r}_{n-k},
\end{align}
and
\be
\left\|\sum\limits_{k=1}^n e^{-i(k-1)\tau{\mathbf{T}}}\Lambda_\tau\hat{r}_{n-k}\right\|_{H^1}\lesssim
\sum\limits_{k=1}^n\frac{\tau^3}{\eps^4}(n-k)\tau\|\Phi_0\|_{L^\infty([0,2\pi\eps^2];(H^1)^2)}\lesssim\tau^2.
\ee
It remains to estimate the exponential sum term. By the decomposition of $\mathbf{T}$ in \eqref{eq:Tdec}, we have
\be
\sum\limits_{k=1}^ne^{-i(k-1)\tau{\mathbf{T}}}\Lambda_\tau e^{-it_{n-k}{\mathbf{T}}} \Phi_0=
\sum\limits_{k=1}^ne^{-i(k-1)\tau{\mathbf{T}_1}}\Lambda_\tau e^{-it_{n-k}{\mathbf{T}_1}} \Phi_0+E_{\rm Res}\Phi_0.
\ee
Since $e^{-is{\mathbf{T}}}=e^{-is{\mathbf{T}_1}}+O(s\mathbf{T}_2)=e^{-is{\mathbf{T}_1}}+O(s\eps^2)$, the residual term can be bounded as
\be
\|E_{\rm Res}\Phi_0\|_{H^1}\lesssim  \sum\limits_{k=1}^nn\tau\eps^2\frac{\tau^3}{\eps^4}\|\Phi_0\|_{H^3}\lesssim\eps^2\tau^2.
\ee
Thus, we get from \eqref{eq:global} that
\be
\left\|{\cal S}_{\tau}^{n}\Phi_0-{\cal S}_{e,\tau}^{n}\Phi_0-\sum\limits_{k=1}^ne^{-i(k-1)\tau{\mathbf{T}_1}}\Lambda_\tau e^{-it_{n-k}{\mathbf{T}_1}} \Phi_0\right\|_{H^1}
\lesssim \tau^2,\quad 1\leq n\leq N=\frac{2\pi\eps^2}{\tau}.
\ee
A direct application of the bounds for $\Lambda_\tau$ in \eqref{eq:lambdabd} leads to
\be
\left\|\sum\limits_{k=1}^ne^{-i(k-1)\tau{\mathbf{T}_1}}\Lambda_\tau e^{-it_{n-k}{\mathbf{T}_1}} \Phi_0\right\|_{H^1}\lesssim\sum\limits_{k=1}^n\frac{\tau^3}{\eps^4}\lesssim\frac{\tau^2}{\eps^2},
\ee
which implies that
\be\label{eq:error1}
\left\|{\cal S}_{\tau}^{n}\Phi_0-{\cal S}_{e,\tau}^{n}\Phi_0\right\|_{H^1}\lesssim \frac{\tau^2}{\eps^2},\quad 1\leq n\leq N=\frac{2\pi\eps^2}{\tau}.
\ee
This result \eqref{eq:error1} will lead to Theorem \ref{thm:tsfp1} and we are going to prove Theorem \ref{thm:tsfp} for better convergence results.

Next, we want to show that for $n=N$, i.e. in one period, the error bounds above can be refined.
The key is to estimate $\sum\limits_{k=1}^Ne^{-i(k-1)\tau{\mathbf{T}_1}}\Lambda_\tau e^{-it_{n-k}{\mathbf{T}_1}} \Phi_0$ in an appropriate  way. Recalling \eqref{eq:Fmid} and
 using $\tau=\frac{2\pi\eps^2}{N}$, we
can write
\begin{align}
&\sum\limits_{k=1}^Ne^{-i(k-1)\tau{\mathbf{T}_1}}\Lambda_\tau e^{-it_{N-k}{\mathbf{T}_1}} \Phi_0=\sum\limits_{k=0}^{N-1}e^{-i(N-k-1)\tau{\mathbf{T}_1}}\Lambda_\tau e^{-it_k{\mathbf{T}_1}} \Phi_0\nonumber\\
&=e^{-i(N-1)\tau\mathbf{T}_1}\left(
\tau\sum\limits_{k=0}^{N-1}{\cal F}_{(k+1/2)\tau}\left(\Phi_0\right)-\int_{0}^{t_N}{\cal F}_s(\Phi_0)\,ds\right),
\end{align}
where the error is contained in the term $\tau\sum\limits_{k=0}^{N-1}{\cal F}_{(k+1/2)\tau}\left(\Phi_0\right)-\int_{0}^{t_N}{\cal F}_s(\Phi_0)\,ds$, a midpoint rule approximation for an
integral of a periodic function over one period. Analogous  to the NLSE case \cite{florian}, such error can be refined.
For a general smooth periodic function $f(s)$ with period $P$, we have (cf. \cite{florian}),
\be
\left\|\frac{P}{N}\sum\limits_{k=0}^{N-1}f\left(\frac{k+1/2}{N}P\right)-\int_0^Pf(s)ds\right\|\lesssim P\left(\frac{P}{2N\pi}\right)^{m_1}\|\partial_s^{m_1}f\|,\quad m_1\ge1.
\ee
Since ${\cal F}_s(\Phi_0)$ is $2\pi\eps^2$ periodic and smooth in $s$, recalling the regularity assumptions (E) and (F), together with the fact that $\mathbf{T}_1$ is bounded from $H^m_p\to H_p^m$ for any $m\in{\mathbb N}$, we find that
\be\label{eq:cancel}
\left\|\sum\limits_{k=1}^Ne^{-i(k-1)\tau{\mathbf{T}_1}}\Lambda_\tau e^{-it_{N-k}{\mathbf{T}_1}} \Phi_0\right\|_{H^1}
\lesssim \eps^2 \tau^{m_\ast},\quad m_\ast\ge1.
\ee
Thus we obtain the refined global error at $n=N$, i.e.
\be\label{eq:errorref}
\left\|{\cal S}_{\tau}^{N}\Phi_0-{\cal S}_{e,\tau}^{N}\Phi_0\right\|_{H^1}\lesssim \tau^2+\eps^2N^{-m_\ast},\quad m_\ast\ge1.
\ee

{\bf Step 3} (bounds on the global error). We are ready to estimate the global error at arbitrary $t_n\leq T$, based on estimates \eqref{eq:error1} and \eqref{eq:errorref}. Let $t_n=2k\pi\eps^2+m\tau$, where $k\ge1$ and $0\leq m\leq N-1$.
Denote the flow
\be
\widetilde{\cal S}_{N}\Psi={\cal S}_{\tau}^N\Psi, \quad \widetilde{\cal S}_{e,N}\Psi={\cal S}_{e,\tau}^N\Psi,
\ee
and it is easy to verify  the stability of $\widetilde{\cal S}_N$ analogous to \eqref{eq:stab}. By the telescopic identity, we have
\begin{align*}
{\cal S}_\tau^n\left(\Phi_0\right)-{\cal S}_{e,\tau}^n(\Phi_0)=&
{\cal S}_\tau^{m}\circ\left(\widetilde{\cal S}_{N}^{k}\left(\Phi_0\right)-\widetilde{\cal S}_{e,N}^{k}\left(\Phi_0\right)\right)-({\cal S}_\tau^m-{\cal S}_{e,\tau}^m)\circ\widetilde{\cal S}_{e,N}^{k}(\Phi_0).
\end{align*}
Applying  the regularity assumption (F) and the error estimate \eqref{eq:error1}, we have
\be
\left\|({\cal S}_\tau^m-{\cal S}_{e,\tau}^m)\circ\widetilde{\cal S}_{e,N}^{k}(\Phi_0)\right\|_{H^1}\lesssim\frac{\tau^2}{\eps^2},\quad 0\leq m\leq N-1.
\ee
Using telescopic identity,  the stability of ${\cal S}_{\tau}$ and estimates \eqref{eq:errorref}, we get
\begin{align*}
\left\|{\cal S}_\tau^{m}\circ\left(\widetilde{\cal S}_{N}^{k}\left(\Phi_0\right)-\widetilde{\cal S}_{e,N}^{k}\left(\Phi_0\right)\right)\right\|_{H^1}
\leq &e^{m\tau C_1}\sum\limits_{j=1}^k\left\|\widetilde{\cal S}_{N}^{j-1}\circ\left(\widetilde{\cal S}_N\circ\widetilde{\cal S}_{e,N}^{k-j}\right)\Phi_0-\widetilde{\cal S}_{N}^{j-1}\circ\left(\widetilde{\cal S}_{e,N}^{k-j+1}\right)\Phi_0\right\|_{H^1}\\
\lesssim&\sum_{j=1}^k(\tau^2+\eps^2h^{m_\ast})\lesssim \frac{T}{2\pi\eps^2}(\tau^2+\eps^2N^{-m_\ast})\lesssim \frac{\tau^2}{\eps^2}+N^{-m_\ast}.
\end{align*}
Thus we have proved the error estimates  \eqref{eq:tssptd1} for the semi-discretization  \eqref{eq:tssptd}.

{\bf Part 2} (convergence of the full discretization). Noticing that
\be
I_M(\Phi^n)(x)-\Phi^{[n]}(x)=I_M(\Phi^n)(x)-P_M(\Phi^{[n]}(x))+P_M(\Phi^{[n]}(x))-\Phi^{[n]}(x),
\ee
we find from the regularity assumption that
\be
\|I_M(\Phi^n)(\cdot)-\Phi^{[n]}(\cdot)\|_{H^s}\leq\|I_M(\Phi^n)(\cdot)-P_M(\Phi^{[n]}(\cdot))\|_{H^s}+\tilde{C}_1h^{m_0-s},\quad s=0,1,
\ee
where $\tilde{C}_1$ is a constant independent of $h$, $n$, $\tau$ and $\eps$.
Hereafter, all the constants used in the inequalities
are independent of $h$, $n$, $\tau$ and $\eps$. We also have error bounds \eqref{eq:tssptd}, i.e.
\be\label{eq:semibd}
\|\Phi^{[n]}(\cdot)-\Phi(t_{n},\cdot)\|_{H^1}\leq \tilde{C}_2(\tau^2/\eps^2+N^{-m_\ast}).
\ee
It suffices to study the error $\bee^n(x)\in Z_M\times Z_M$ given as
\be
\bee^n(x)=I_M(\Phi^n)(x)-P_M(\Phi^{[n]})(x),\quad 0\leq n\leq T/\tau.
\ee

We shall prove \eqref{eq:errortsfp} by mathematical induction, i.e. for $0\le n\leq\frac{T}{\tau}$,
\be\label{eq:err1}
\|I_M(\Phi^n)(x)-\Phi(t_n,x)\|_{H^s}\leq C_{T}(\tau^2/\eps^2+N^{-m_\ast})+C_Hh^{m_0-s},\quad \|\Phi^n\|_{l^\infty}\leq M_0+1, \quad s=0,1,
\ee
where $C_T$ and $C_H$ (independent of $h$, $n$, $\eps$ and $\tau$) are constants to be determined later.

It is easy to check that when $n=0$, we have $\|\bee^0(x)\|_{H^s}\leq \tilde{C}_3h^{m_0-s}$ ($s=0,1$) and the estimates \eqref{eq:errortsfp} hold if $C_H\ge\max\{\tilde{C}_1,\tilde{C}_3\}$ and $h$ is small enough.

Assume that for $0\leq n\leq m\leq\frac{T}{\tau}-1$, the error estimates \eqref{eq:err1} hold. For  $n=m+1$, we have for $\Phi^{m+1}$ and $\Phi^{[m+1]}$
\beas
&&I_M(\Phi^{(1)})=e^{-i\tau{\mathbf{T}}/2}I_M(\Phi^{m}),
\quad I_M(\Phi^{m+1})=e^{-i\tau{\mathbf{T}}/2}I_M(\Phi^{(2)}),\\
&&I_M(\Phi^{(2)})=I_M(e^{-i\tau(V(x_j)I_2-A_1(x_j)\sigma_1+\lambda_2|\Phi^{(1)}_j|)}\Phi^{(1)}_j),\\
&&P_M(\Phi^{<1>})=e^{-i\tau{\mathbf{T}}/2}I_M(\Phi^{[m]}),\quad P_M(\Phi^{[m+1]})=e^{-i\tau{\mathbf{T}}/2}P_M(\Phi^{<2>}),\\
&&P_M(\Phi^{<2>})=P_M(e^{-i\tau(V(x)I_2-A_1(x)\sigma_1+\lambda_2|\Phi^{<1>}|^2)}\Phi^{<1>}).
\eeas
As $e^{-i\tau\mathbf{T}}$ preserves $H^s$ norm, we get
\bea
\|\bee^m(\cdot)\|_{H^s}=\|I_M(\Phi^{(1)})-P_M(\Phi^{<1>})\|_{H^s},\ \ \|\bee^{m+1}(\cdot)\|_{H^s}=\|I_M(\Phi^{(2)})-P_M(\Phi^{<2>})\|_{H^s},\ \ s=0,1. \quad
\eea
On the other hand, we have
\begin{equation*}
I_M(\Phi^{(2)})- P_M(\Phi^{<2>})=I_M(e^{-i\tau(V(x_j)I_2-A_1(x_j)\sigma_1+\lambda_2|\Phi^{(1)}_j|)}
\Phi^{(1)}_j)-P_M(e^{-i\tau(V(x)I_2-A_1(x)\sigma_1+\lambda_2|\Phi^{<1>}|^2)}\Phi^{<1>}),
\end{equation*}
which together with $\Phi^{<1>}\in H_p^{m_0}$ implies
\begin{align}
\|I_M(\Phi^{(2)})- P_M(\Phi^{<2>})\|_{H^s}\leq \tilde{C}_3h^{m_0-s}+\|W(x)\|_{H^s},
\end{align}
where
\be
W(x):=I_M(e^{-i\tau(V(x_j)I_2-A_1(x_j)\sigma_1+
\lambda_2|\Phi^{(1)}_j|)}\Phi^{(1)}_j)-I_M(e^{-i\tau(V(x)I_2-A_1(x)\sigma_1+\lambda_2|\Phi^{<1>}|^2)}\Phi^{<1>}).
\ee
As shown in \cite{BC1,BC2,BC3,BCZ}, $W(x)$ can be estimated through finite difference approximation as
\beas
&&\|W(x)\|_{L^2}\leq C\tau\left\|\Phi^{(1)}_j-\Phi^{<1>}(x_j)\right\|_{l^2}\leq \tilde{C}_4\tau\left(\|\bee^{m}(\cdot)\|_{L^2}+h^{m_0}\right),\\
&&\|W(x)\|_{H^1}\leq C\tau\left(\left\|\Phi^{(1)}_j-\Phi^{<1>}(x_j)\right\|_{l^2}+
\left\|\delta_{x}^+(\Phi^{(1)}_j-\Phi^{<1>}(x_j))\right\|_{l^2}\right)\leq \tilde{C}_5\tau\left(\|\bee^{m}(\cdot)\|_{H^1}+h^{m_0-1}\right),\quad
\eeas
where $\delta_x^+\Phi_j=\frac{\Phi_{j+1}-\Phi_j}{h}$ is the forward finite difference operator. The key point is that $\|\partial_x(I_M\Psi_j)\|_{L^2}\sim \|\delta_x^+\Psi_j\|_{l^2}$.
Thus, we have
\be
\|\bee^{m+1}(\cdot)\|_{H^s}\leq \max\{\tilde{C}_4,\tilde{C}_5\}\tau\left(\|\bee^{m}(\cdot)\|_{H^s}+h^{m_0-s}\right),\quad s=0,1.
\ee
Indeed, it is true for all $n\leq m$,
\be
\|\bee^{n+1}(\cdot)\|_{H^s}\leq \max\{\tilde{C}_4,\tilde{C}_5\}\tau\left(\|\bee^{n}(\cdot)\|_{H^s}+h^{m_0-s}\right),\quad s=0,1.
\ee
Using discrete Gronwal inequality, we get
\be
\|\bee^{n+1}(\cdot)\|_{H^s}\leq \tilde{C}_6 h^{m_0-s},\quad n\leq m\leq \frac{T}{\tau}-1.
\ee
Thus \eqref{eq:err1} holds true for $n=m+1$ if we choose $C_T=\tilde{C}_2$, $C_H=\max\{\tilde{C}_1,\tilde{C}_3,\tilde{C}_6\}$ and use the discrete Sobolev inequality with
sufficiently small $h$ and $\tau$.
This completes the induction and Theorem \ref{thm:tsfp} holds.

\Acknowledgements{This work was partially supported by the Ministry
of Education of Singapore grant
R-146-000-196-112 (W. Bao) and the Natural Science Foundation
of China Grant 91430103 (Y. Cai).
Part of this work was done when the authors were visiting
the Institute for Mathematical Sciences at the National University of Singapore in 2015.}


\end{document}